\documentclass[10pt]{amsart}
\usepackage{}
\usepackage{latexsym,amsfonts,amssymb,mathrsfs,graphicx,subfigure,stmaryrd}
\usepackage{color}
\usepackage{accents}

\evensidemargin=\oddsidemargin
\newdimen\dummy
\dummy=\oddsidemargin
\newcommand{\Lr}[1]{\left(#1\right)}

\def\al{\alpha}
\def\om{\Omega}

\def\pa{\partial}

\def\ha{a'}

\def\whw{\widehat{W}}
\def\wtw{\widetilde{W}}
\def\wta{\widetilde{A}}
\def\wts{\widetilde{S}}

\def\E{\mathcal E}

\def\T{\mathcal T}
\def\H{\mathcal {H}}

\def\na{\nabla}

\def\P{\mathcal {P}}
\def\R{\mathbb R}

\def\V{\mathcal {V}}

\newtheorem{theorem}{Theorem}[section]
\newtheorem{lemma}[theorem]{Lemma}
\newtheorem{corollary}[theorem]{Corollary}

\newtheorem{example}[theorem]{Example}

\theoremstyle{definition} \numberwithin{equation}{section}

\begin {document}
\title [FETI-DP of DG for Multiscale Problems]{A FETI-DP Preconditioner of Discontinuous Galerkin Method For Multiscale Problems in High contrast Media}

\author[R. Du]{Rui Du}
\address{Department of Informatics, University of Bergen, N-5020 Bergen, Norway}\email{Rui.Du@ii.uib.no}

\author[Y.-F. Ma]{Yunfei Ma}
\address{Institute of
Computational Mathematics and Scientific/Engineering Computing,
Academy of Mathematics and Systems Science, Chinese Academy of Sciences, 100190 Beijing, P.R. China}
\email{mayf@lsec.cc.ac.cn}

\author[T. Rahman]{Talal Rahman}
\address{Department of Computing, Mathematics and Physics, Faculty of Engineering and Business Administration, Bergen University College, N-5020 Bergen, Norway}
\email{Talal.Rahman@hib.no}

\author[X.-J. Xu]{Xuejun Xu}
\address{Institute of
Computational Mathematics and Scientific/Engineering Computing,
Academy of Mathematics and Systems Science, Chinese Academy of Sciences, Beijing, 100190, P.R. China}
\email{xxj@lsec.cc.ac.cn}

\thanks{The authors are thankful to M. Dryja for many fruitful discussions on the topics of this paper.}
%    General info
\subjclass[]{}
%\date{\today}

\begin {abstract}

 In this paper we consider the second order elliptic partial differential equations with highly varying (heterogeneous) coefficients on a two-dimensional region. The problems are discretized by a composite finite element (FE) and discontinuous Galerkin (DG) Method. The fine grids are in general nonmatching across the subdomain boundaries, and the subdomain partitioning does not need to resolve the jumps in the coefficient. A FETI-DP preconditioner is proposed and analyzed to solve the resulting linear system. Numerical results are presented to support our theory.

\end {abstract}
\maketitle

\section{Introduction}

We consider the following problem: Find $u^* \in H_0^1(\om)$ such that
\begin {equation}\label {Dip}
\begin {aligned}
a(u^*, v)=(f, v)\qquad\text{for all}\qquad v\in H_0^1(\om),
\end {aligned}
\end {equation}
where
\begin {align*}
a(u,v):= \int_{\om} \al(x)\na u\cdot\na v dx \quad\text{and}\quad (f,v):=\int_{\om} fv dx,
\end {align*}
where $\om\subset\R^2$ is a bounded polygonal domain. We assume that $\al(x) \geq \al_0 >0$ and $\al(x)\in L^\infty(\om)$ which may be discontinuous, while $f(x)\in L^2(\om)$. The representative examples of the problem~\eqref {Dip} are subsurface flows in heterogeneous media~\cite {LJ06, LJ08} where the heterogeneity varies over a wide range of scales. The aim of this paper is to design and analyze a FETI-DP method for solving such problems based on a composite FE/DG discretization.

Instead of using the full DG method over the whole domain, the composite FE/DG method employs conforming FE methods inside the subdomains, while applies a DG discretization only on the subdomain interfaces to deal with the nonmatching meshes across the interfaces; see~\cite {ABCM02, DDG03, DGS07, DGS13, EGLMS13}. The local bilinear forms of the discrete problem are composed of three symmetric terms: the one associated with the energy, the one ensuring consistency and symmetry, and the interior penalty term~\cite {St98, Ri08} to handle the nonconforming FE spaces across the interfaces; see cf.~\eqref {bi1}-~\eqref {bi4}.
Such discretization allows for nonmatching grids which provides greater flexibility in the choice of mesh partitioning and memory storage. This may be useful particularly when the coefficient varies roughly in one subdomain and mildly in the others.

FETI-DP methods, as well as FETI~\cite {FR91, FMR94} and BDDC~\cite {Do03, MD03}, have been well established as a class of nonoverlapping domain decomposition methods for solving large-scale linear systems. These methods have been used widely for standard continuous FE discretization, and verified to be successful both theoretically and numerically; see ~\cite {TW05} and references therein. FETI-DP method was firstly introduced in~\cite {FL01} following by a theoretical analysis provided in~\cite {MT01}. In FETI-DP algorithms, we need a relatively small number of continuity constraints across the interface in each iteration step. The continuity of the solution across the subdomain interfaces is enforced by Lagrange multipliers, while the continuity at the subdomain vertices is enforced directly by assigning unique values. The methods were further improved in~\cite {FLP00, KWD02, TW05} to use the continuity constraints on the averages across the edges on subdomain interfaces. The FETI-DP methods have been developed more recently, and possess several advantages over the one-level FETI method; see cf.~\cite {TW05}.

The FETI-DP method was firstly considered for composite FE/DG discretization in~\cite {DGS13}. We will follow the same framework as described therein. In~\cite {DGS13}, the discontinuities of the coefficients are assumed to occur only across the subdomain interfaces. The main purpose of this paper is to extend the methodology to the case where the coefficients are allowed to have large jumps not only across but also along the subdomain interfaces and in the interior of the subdomains. We recall that such problems were investigated in the context of FETI methods in~\cite {PS08, PS11}.

In this paper, we will use the same DG bilinear form as in~\cite {EGLMS13}, construct our FETI-DP preconditioner as in~\cite {DGS13}, and define the components of the scaling matrix as proposed in~\cite {PS08}. For the theoretical aspect, we employ the cut off technique and the generalized discrete Sobolev type inequality, cf.~\cite {EGLMS13}, as well as the standard estimates of the edge and vertex functions, cf.~\cite {TW05}. It will be proved that the convergence of the FETI-DP method only weakly depends on the jump of coefficients, i.e., linearly depends on the contrast of the coefficients in the boundary layer. Here we define the boundary layer as the union of fine triangles that touch the subdomain boundaries. We also show that the condition number estimate of the proposed method is quadratic dependence on $H/h$ where $H$ is the subdomain diameter and $h$ is the fine mesh size. This quadratic dependence on $H/h$ can be relaxed to a weaker dependence of $H/h(1+\log H/h)^2$ under stronger assumptions on the coefficients in the interior of the subdomains.

The remaining part of this paper is organized as follows. In Section 2, we introduce the composite FE/DG formulation of problem~\eqref {Dip}. The FETI-DP method is presented in Section 3. The main results of the paper are given in Section 4 about the analysis of the condition number estimate. Numerical results are provided in Section 5 to confirm the theoretical analysis. In the last section we summarize our findings and discuss certain extensions.

Throughout this paper we denote a Sobolev space
of order $k$ by the standard notation $H^k(\om)$ with norm given by $\|\cdot\|_{H^k(\om)}$; see e.g.,~\cite {Af03} for exact definition. For $k=0$ we
use $L^2(\om)$ instead of $H^0(\om)$ and write the norm as $\|\cdot\|_{L^2(\om)}$. In addition, $A \simeq B$ stands for $C_1 B \leq A \leq C_2 B$ with positive constants $C_1$ and $C_2$ depending only on the shape regularity of the meshes.

\section{DG Discretization} In this section we present the DG formulations of problem~\eqref {Dip} that will be studied here.

Let the domain $\overline{\om} = \cup_{i=1}^N \overline{\om}_i$ and $\om_i$ be disjoint shape regular polygonal subdomains of diameters $H_i$. Denote the subdomain boundaries by $\pa\om_i$. For each $\om_i$, we introduce a shape regular triangulation $\T_h(\om_i)$ with the mesh size $h_i$. Note that the resulting triangulation of $\om$ is in general nonmatching across $\pa\om_i$.

%Let $\T_h(\om)$ be a fine triangulation on $\om$ with the mesh size $h$. It is constructed from a refinement of subdomains $\om_i$. We assume that the fine triangulation $\T_h(\om)$ is matching across $\pa\om_i$.

We assume that the substructures $\{\om_i\}_{i=1}^N$ form a geometrically conforming partition of $\om$, i.e., the intersection $\pa\om_i \cap \pa\om_j (i\neq j)$ is either empty, or a common vertex or edge of $\om_i$ and $\om_j$. Let us denote the common edge $\bar{E}_{ij} = \bar{E}_{ji} := \pa\om_i \cap \pa\om_j$. Although $E_{ij}$ and $E_{ji}$ are geometrically the same object, we will treat them separately since we consider different triangulations on $\bar{E}_{ij}\subset\pa\om_i$ and on $\bar{E}_{ji}\subset\pa\om_j$, with the mesh size of $h_i$ and $h_j$, respectively. In the text below, we use $E_{ijh}$ and $E_{jih}$ to denote the set of nodal points of the triangulation on $E_{ij}$ and $E_{ji}$ with mesh sizes $h_i$ and $h_j$, respectively, and $\bar{E}_{ijh}$ and $\bar{E}_{jih}$ when the endpoints are included. Moreover, the two triangulations $\T_h(\om_i)$ and $\T_h(\om_j)$ can be merged to obtain a finer but the same mesh on $\bar{E}_{ij}$ and $\bar{E}_{ji}$.

We also denote $E_{i\pa}:= \pa\om_i \cap \pa\om$ when there is an intersection between $\pa\om_i$ and the global boundary $\pa\om$. Let us denote by $\E_i^0$ the set of indices to refer to the edges $E_{ij}$, i.e., $j$ of $\om_j$ which has a common edge $E_{ji}$ with $\om_i$, and by $\E_i^\pa$ the set of indices to refer to the edges $E_{i\pa}$. The set of indices of all edges of $\om_i$ is denoted by $\E_i := \E_i^0\cup\E_i^\pa$.

For simplicity, we assume that the coefficient $\al(x)\geq \alpha_0 = 1$, which can be fulfilled by scaling~\eqref {Dip} with $1/\min_x \al(x)$. Without loss of generality again, we assume that $\al(x)$ is constant over each fine triangle. The analysis will depend on the coefficient in a boundary layer near subdomain boundaries. For each subdomain $\om_i$, we define the boundary layer $\om_i^h$ by
\begin {align*}
\overline{\om_i^h}:=\bigcup\{\bar{\tau}: \tau\in \overline{\om}_i, \text{dist}(\tau, \pa\om_i) \leq h_i\},
\end {align*}
i.e., the union of fine triangles in $\mathcal {T}_h(\om_i)$ that touch the boundary $\partial\om_i$. Furthermore, we set
\begin {equation}\label {coef}
\begin {aligned}
\underline{\al}_i := \inf_{x\in\om_i^h} \al(x)\quad\text{and}\quad
\overline{\al}_i := \sup_{x\in\om_i^h} \al(x).
\end {aligned}
\end {equation}
Let $\al_i(x)$ be $\al(x)$ restricted to $\overline{\om}_i$. We define the harmonic averages along the edges $E_{ij}$ as follows:
\begin {equation}\label {haraver}
\al_{ij}(x) = \dfrac{2 \al_i(x)\al_j(x)}{\al_i(x)+\al_j(x)}\quad\text{and}\quad h_{ij}= \dfrac{2 h_i h_j}{h_i+h_j}.
\end {equation}
Note that the functions $\al_{ij}(x)$ and $h_{ij}$ are piecewise constant over the edge $E_{ij}$ on the mesh that is obtained by merging the partitions $\T_h(\om_i)$ and $\T_h(\om_j)$ along this common edge $E_{ij}$. It is easy to check that
\begin {equation}\label {ineqaver}
\min(\underline{\al}_i , \underline{\al}_j)\leq \al_{ij} \leq 2\min(\underline{\al}_i, \underline{\al}_j)\quad\text{and}\quad \min(h_i, h_j)\leq h_{ij} \leq 2\min(h_i, h_j).
\end {equation}

Let $V_h(\om_i)$ be the standard finite element space of continuous piecewise linear functions in $\om_i$. Define
\begin {equation}\label {space}
\begin {aligned}
V_h(\om)=\prod_{i=1}^N V_h(\om_i)\equiv V_h(\om_1)\times V_h(\om_2)\times\cdots\times V_h(\om_N),
\end {aligned}
\end {equation}
and represent functions $u\in V_h(\om)$ as $u = \{u_i\}_{i=1}^N$ with $u_i\in V_h(\om_i)$.
We do not assume that functions in $V_h(\om_i)$ vanish on $\pa\om_i\cap\pa\om$.

The discrete problem obtained by the DG method is of the form: Find $u_h^*=\{u_{h,i}^*\}_{i=1}^N\in V_h(\om)$ with $u_{h,i}^*\in V_h(\om_i)$ such that
\begin {equation}\label {disp}
\begin {aligned}
a_h(u_h^*,v)=(f,v)\qquad\text{for all}\qquad v=\{v_i\}_{i=1}^N\in V_h(\om),
\end {aligned}
\end {equation}
where
\begin {align*}
a_h(u,v):=\sum_{i=1}^N a'_i(u,v)\quad\text{and}\quad (f,v):=\sum_{i=1}^N\int_{\om_i}fv_i dx.
\end {align*}
Here each local bilinear form $a'_i(\cdot,\cdot)$ is given as the sum of three symmetric bilinear forms:
\begin {equation}\label {bi1}
\begin {aligned}
a'_i(u,v):=a_i(u,v)+s_i(u,v)+p_i(u,v),
\end {aligned}
\end {equation}
where
\begin {equation}\label {bi2}
\begin {aligned}
%a_i(u,v):=\sum_{\tau\subset\om_i}\al_\tau\int_{\tau}\na u_i\cdot\na v_idx,
a_i(u,v):=\int_{\om_i}\al_i(x)\na u_i\cdot\na v_idx,
\end {aligned}
\end {equation}
\begin {equation}\label {bi3}
\begin {aligned}
%s_i(u,v):=\sum_{j\in\E_i}\sum_{e\subset E_{ij}}\frac{\al_e}{l_{ij}}\int_{e}\Lr{\frac{\pa u_i}{\pa n}(v_j-v_i)+\frac{\pa v_i}{\pa n}(u_j-u_i)}ds,
s_i(u,v):=\sum_{j\in\E_i}\frac{1}{l_{ij}}\int_{E_{ij}}\al_{ij}(x)\Lr{\frac{\pa u_i}{\pa n_i}(v_j-v_i)+\frac{\pa v_i}{\pa n_i}(u_j-u_i)}ds,
\end {aligned}
\end {equation}
and
\begin {equation}\label {bi4}
\begin {aligned}
%p_i(u,v):=\sum_{j\in\E_i}\sum_{e\subset E_{ij}}\frac{\delta\al_e}{l_{ij}h_{ij}}\int_{e}(u_j-u_i)(v_j-v_i)ds,
p_i(u,v):=\sum_{j\in\E_i}\frac{1}{l_{ij}}\frac{\delta}{h_{ij}}\int_{E_{ij}}\al_{ij}(x)(u_j-u_i)(v_j-v_i)ds.
\end {aligned}
\end {equation}
Here $\frac{\pa}{\pa n_i}$ denotes the outward normal derivative on $\pa\om_i$, and $\delta$ is a positive penalty parameter. When $j\in\E_i^0$, we set $l_{ij} = 2$, and let $\al_{ij}$ and $h_{ij}$ be defined in~\eqref {haraver}. When $j\in\E_i^\pa$, we set $l_{i\pa}=1, u_\pa=0, v_\pa=0$, and define $\al_{i\pa} = \al_i$ and $h_{i\pa}=h_i$.

We introduce the bilinear form
\begin {equation}\label {nbi1}
\begin {aligned}
d_h(u,v):=\sum_{i=1}^N d_i(u,v)
\end {aligned}
\end {equation}
with
\begin {equation}\label {nbi2}
\begin {aligned}
d_i(u,v):=a_i(u,v)+p_i(u,v).
\end {aligned}
\end {equation}
It is easy to check that $d_h(\cdot,\cdot)$ is symmetric and positive definite, which can induce a broken norm in $V_h(\om)$ by
\begin {align*}
%||u||_h^2:=d_h(u,u)=\sum_{i=1}^N\Lr{\sum_{\tau\subset\om_i}\alpha_\tau||\na u_i||_{L^2(\tau)}^2+\sum_{j\in\E_i}\sum_{e\subset E_{ij}}\frac{\delta\al_e}{l_{ij}h_{ij}}||u_i-u_j||_{L^2(e)}^2},
||u||_h^2:=d_h(u,u)=\sum_{i=1}^N\Lr{||\al_i^{1/2}\na u_i||_{L^2(\om_i)}^2+\sum_{j\in\E_i}\frac{1}{l_{ij}}\frac{\delta}{h_{ij}}||\al_{ij}^{1/2}(u_i-u_j)||_{L^2(E_{ij})}^2}
\end {align*}
for any $u=\{u_i\}_{i=1}^N\in V_h(\om)$.

The next lemma characterizes the equivalence between the bilinear forms $a_h(\cdot,\cdot)$ and $d_h(\cdot,\cdot)$. This equivalence implies the existence and uniqueness of the solution to the discrete problem~\eqref {disp}, and also allows us to use the bilinear form $d_h(\cdot,\cdot)$ instead of $a_h(\cdot,\cdot)$ for preconditioning.
\begin {lemma}\label {equibi}
There exists $\delta_0>0$ such that for $\delta\geq \delta_0$ and for all $u\in V_h(\om)$, we have
\begin {equation}\label {equibi1}
\begin {aligned}
\gamma_0 d_i(u,u)\leq \ha_i(u,u)\leq \gamma_1 d_i(u,u)\quad\text{for all}\quad i=1,\cdots,N,
\end {aligned}
\end {equation}
and
\begin {equation}\label {equibi2}
\begin {aligned}
\gamma_0 d_h(u,u)\leq a_h(u,u)\leq \gamma_1 d_h(u,u),
\end {aligned}
\end {equation}
where $\gamma_0$ and $\gamma_1$ are positive constants independent of $h_i$, $H_i$, $\al_i(x)$, and $u$. For the proof we refer to Lemma 2.1 of~\cite {EGLMS13}.
\end {lemma}

\section {FETI-DP Preconditioner for the Schur Complement Systems}
\label {FETI-DP for DG}

In this section, we will give the formulation of our FETI-DP method using the framework introduced in~\cite {TW05, DGS13}.

\subsection {Schur Complement Systems and Discrete Harmonic Extensions}

Firstly, we borrow the notations from~\cite {DGS13}. Let
\begin {align*}
\om'_i:=\overline{\om}_i\bigcup\{\cup_{j\in\E_i^0}\bar{E}_{ji}\},
\end {align*}
i.e., the union of $\overline{\om}_i$ and the $\bar{E}_{ji}\subset\pa\om_j$ with $j\in\E_i^0$, and let
\begin {align*}
\Gamma_i:=\overline{\pa\om_i\backslash\pa\om},\quad
\Gamma'_i:=\Gamma_i\bigcup\{\cup_{j\in\E_i^0}\bar{E}_{ji}\}, \quad\text{and}\quad I_i:=\om'_i\backslash\Gamma'_i.
\end {align*}
Then we set
\begin {align}\label {notation}
\Gamma:=\bigcup_{i=1}^N\Gamma_i, \quad \Gamma':=\prod_{i=1}^N\Gamma'_i,\quad\text{and}\quad I:=\prod_{i=1}^N I_i.
\end {align}

% A function $u_i\in W_i(\om'_i)$ can be represented as
%\begin {align}\label {IBdecomp}
%u_i = \{(u_i)_i,\{(u_i)_j\}_{j\in\E_i^0}\},
%\end {align}
%where $(u_i)_i$ is $u_i$ restricted to ${\overline{\om}_i}$ and $(u_i)_j$ is $u_i$ restricted to ${\bar{E}_{ji}}$ with $j\in \E_i^0$.

We introduce $W_i(\om'_i)$ as the FE space of functions defined on the nodal values of $\om'_i$. That is,
\begin {equation}\label {wspacelocal1}
\begin {aligned}
W_i(\om'_i) = W_i(\overline{\om}_i)\times \prod_{j\in\E_i^0} W_i(\bar{E}_{ji}),
\end {aligned}
\end {equation}
where $W_i(\bar{E}_{ji})$ is the trace of the space $V_h(\om_j)$ on $\bar{E}_{ji}\subset\pa\om_j$ with $j\in\E_i^0$. In the following, we use the same notation to denote both FE functions and their vector representations. The local bilinear form $a'_i(\cdot,\cdot)$ in~\eqref {bi1} is defined over $W_i(\om'_i)\times W_i(\om'_i)$, and the associated stiffness matrix is given by
\begin {equation}\label {bistiff}
\langle A'_i u_i, v_i\rangle=a'_i(u_i, v_i)\quad\text{for all}\quad u_i, v_i\in W_i(\om'_i),
\end {equation}
where $\langle \cdot, \cdot\rangle$ denotes the Euclidean inner product associated to the vectors with nodal values. We will decompose $u_i\in W_i(\om'_i)$ as $u_i = (u_{i,I}, u_{i,\Gamma'})$, where $u_{i,I}$ represents values of $u_i$ at interior nodal points on $I_i$ and $u_{i,\Gamma'}$ at the nodal points on $\Gamma'_i$. Note that for subdomains $\om_i$ which intersect $\pa\om$ by edges, the nodal values of $W_i(\om'_i)$ on $\pa\om_i\backslash\Gamma'_i$ are treated as unknowns and belong to $I_i$. Hence, we can rewrite
\begin {align}\label {wspacelocal2}
W_i(\om'_i) = W_i(I_i)\times W_i(\Gamma'_i),
\end {align}
and the matrix $A'_i$ as
\begin {equation}\label {stiff}
A'_i=\left(
            \begin{array}{cc}
A'_{i,II}& A'_{i,I\Gamma'}\\
A'_{i,\Gamma' I}& A'_{i,\Gamma'\Gamma'}
\end{array}
          \right),
\end {equation}
where the block rows and columns correspond to the nodal points of $I_i$ and $\Gamma'_i$, respectively.

The Schur Complement $S'_i$ of $A'_i$, with respect to the nodal points of $\Gamma'_i$, takes the form
\begin {align}\label {Schurlocal}
S'_i := A'_{i,\Gamma'\Gamma'} - A'_{i,\Gamma' I}(A'_{i,II})^{-1}A'_{i,I\Gamma'}.
\end {align}
Note that $S'_i$ satisfies the energy minimizing property
\begin {align}\label {mini}
\langle S'_i u_{i,\Gamma'}, u_{i,\Gamma'}\rangle=\min a'_i(w_i, w_i)
\end {align}
subject to the condition that $w_i = (w_{i,I}, w_{i,\Gamma'})\in W_i(\om'_i)$ and $w_{i,\Gamma'} = u_{i,\Gamma'}$ on $\Gamma'_i$. The bilinear form $a'_i(\cdot,\cdot)$ is symmetric and nonnegative with respect to $W_i(\om'_i)$, see Lemma~\ref {equibi}. The minimizing function of~\eqref {mini} is called the discrete harmonic extension in the sense of $a'_i(\cdot,\cdot)$, denoted by $\H'_i u_{i,\Gamma'}$, and satisfies
\begin {align}\label {dishar}
a'_i(\H'_i u_{i,\Gamma'}, v_i) = 0\quad\text{for all}\quad v_i\in\accentset{\circ}{W}_i(\om'_i)
\end {align}
with $\H'_i u_{i,\Gamma'} = u_{i,\Gamma'}$ on $\Gamma'_i$. Here $\accentset{\circ}{W}_i(\om'_i)$ is the subspace of $W_i(\om'_i)$ of functions which vanish on $\Gamma'_i$. We also introduce $\H_i u_{i,\Gamma'}\in W_i(\om'_i)$, the standard discrete harmonic extension in the sense of $a_i(\cdot,\cdot)$, which is defined by
\begin {align}\label {disharS}
a_i(\H_i u_{i,\Gamma'}, v_i) = 0\quad\text{for all}\quad v_i\in\accentset{\circ}{W}_i(\om'_i)
\end {align}
with $\H_i u_{i,\Gamma'} = u_{i,\Gamma'}$ on $\Gamma'_i$.

Note that the extensions, $\H_i$ and $\H'_i$, differ from each other in the sense that $\H_i u_{i,\Gamma'}$ at the interior nodes $I_i$ depends only on the nodal values of $u_{i,\Gamma'}$ on $\Gamma_i$ while $\H'_i u_{i,\Gamma'}$ depends on the nodal values of $u_{i,\Gamma'}$ on $\Gamma'_i$. The next lemma shows the equivalence between $\H_i$ and $\H'_i$ in the energy form induced by $d_i(\cdot,\cdot)$. This equivalence will allow us to take advantages of all the discrete Sobolev results known for $\H_i$ discrete harmonic extensions. The fundamental idea of the proof comes from~\cite {DGS07}, and we still include the proof here for completeness.

\begin {lemma}\label {enerequiv}
For any $u_{i,\Gamma'}\in W_i(\Gamma'_i)$, there exists a constant $C > 0$ independent of $h_i, H_i, \al_i(x)$ and $u_{i,\Gamma'}$, such that
\begin {align}\label {dnormequiv1}
d_i(\H_i u_{i,\Gamma'}, \H_i u_{i,\Gamma'})\leq d_i(\H'_i u_{i,\Gamma'}, \H'_i u_{i,\Gamma'})\leq C d_i(\H_i u_{i,\Gamma'}, \H_i u_{i,\Gamma'}).
\end {align}
\end {lemma}
\begin {proof}Here and below, for simplicity of presentation, we omit the subscript $\Gamma'$
and denote $u_{i,\Gamma'}$ by $u_i$ if there is no confusion.

The left-hand inequality of~\eqref {dnormequiv1} follows from the energy minimizing
property of the discrete harmonic extension $\H_i$ in the sense of $a_i(\cdot,\cdot)$,
and the fact that $\H_i u_i = \H'_i u_i = u_i$ on $\Gamma'_i$. Here we remain to prove
the right-hand inequality.

It is easy to verify that $\H_i\H'_i u_i = \H_i u_i$ since the extensions keep the boundary values. Note that we can represent $\H'_i u_i \in W_i(\om'_i)$ as
\begin {align}\label {project}
\H'_i u_i = \H_i \H'_i u_i + \P_i \H'_i u_i,
\end {align}
where $\P_i \H'_i u_i$ is the projection of $\H'_i u_i$ into $\accentset{\circ}{W}_i(\om'_i)$ in the sense of $a_i(\cdot,\cdot)$, i.e., $\P_i \H'_i u_i \in\accentset{\circ}{W}_i(\om'_i)$ and satisfies
\begin {align*}
a_i(\P_i \H'_i u_i, v_i) = a_i(\H'_i u_i, v_i)\quad\text{for all}\quad v_i\in\accentset{\circ}{W}_i(\om'_i).
\end {align*}
Choosing $v_i = \P_i \H'_i u_i$, by Cauchy-Schwarz inequality, we obtain
\begin {align}\label {csineq}
a_i(\P_i \H'_i u_i, \P_i \H'_i u_i) \leq a_i(\H'_i u_i, \H'_i u_i).
\end {align}
Hence,
\begin {equation}\label {enerequiv1}
\begin {aligned}
d_i(\H'_i u_i, \H'_i u_i) &= d_i(\H'_i u_i, \H_i \H'_i u_i) + d_i(\H'_i u_i, \P_i \H'_i u_i)\\
&=d_i(\H'_i u_i, \H_i u_i) + d_i(\H'_i u_i, \P_i \H'_i u_i).
\end {aligned}
\end {equation}

Since the bilinear form $d_i(\cdot,\cdot)$ is symmetric and nonnegative, by Cauchy-Schwarz inequality again, we have
\begin {align}\label {part1}
d_i(\H'_i u_i, \H_i u_i) \leq \epsilon d_i(\H'_i u_i, \H'_i u_i) + \frac{1}{4\epsilon}d_i(\H_i u_i, \H_i u_i)
\end {align}
with arbitrary $\epsilon > 0$.

Since $\P_i \H'_i u_i\in\accentset{\circ}{W}_i(\om'_i)$, using the formulations~\eqref {bi4} and~\eqref {dishar}, we get
\begin {align*}
d_i(\H'_i u_i, \P_i \H'_i u_i) = a_i(\H'_i u_i, \P_i \H'_i u_i),
\end {align*}
and
\begin {align*}
0 = a'_i(\H'_i u_i, \P_i \H'_i u_i) = a_i(\H'_i u_i, \P_i \H'_i u_i) + s_i(\H'_i u_i, \P_i \H'_i u_i),
\end {align*}
which together imply that
\[
d_i(\H'_i u_i, \P_i \H'_i u_i) = - s_i(\H'_i u_i, \P_i \H'_i u_i).
\]
%From the above two equations, we have
%\begin {equation}\label {enerequiv3}
%\begin {aligned}
%d_i(\H'_i u_i, \P_i \H'_i u_i) &= - s_i(\H'_i u_i, \P_i \H'_i u_i)\\
%&= - \sum_{j\in\E_i}\sum_{e\subset E_{ij}}\frac{\al_e}{l_{ij}}\int_{e}\frac{\pa (\P_i \H'_i u_i)}{\pa n}\Lr{(u_i)_j-(u_i)_i}ds
%\end {aligned}
%\end {equation}
We proceed the same lines of Lemma 2.1 in~\cite {EGLMS13}, and finally obtain
\begin {equation}\label {part2}
\begin {aligned}
d_i(\H'_i u_i, \P_i \H'_i u_i) &\leq C \Lr{2\epsilon a_i(\P_i \H'_i u_i, \P_i \H'_i u_i)  + \frac{1}{2\epsilon\delta} p_i(\H_i u_i, \H_i u_i)}\\
&\leq C \Lr{2\epsilon a_i(\H'_i u_i, \H'_i u_i)  + \frac{1}{2\epsilon\delta} p_i(\H_i u_i, \H_i u_i)}\\
&\leq C \Lr{2\epsilon d_i(\H'_i u_i, \H'_i u_i)  + \frac{1}{2\epsilon\delta} d_i(\H_i u_i, \H_i u_i)},
\end {aligned}
\end {equation}
where we have used~\eqref {csineq}.

Combining~\eqref {part1} and~\eqref {part2}, we have
\begin {align*}
d_i(\H'_i u_i, \H'_i u_i) \leq C\Lr{\epsilon d_i(\H'_i u_i, \H'_i u_i) + \frac{1}{4\epsilon} d_i(\H_i u_i, \H_i u_i)}.
\end {align*}
The right-hand side of~\eqref {dnormequiv} follows by choosing a sufficiently small $\epsilon$.
\end {proof}

Lemma~\ref {equibi} and Lemma~\ref {enerequiv} together directly give the following corollary.

\begin {corollary}\label {coro1}
For any $u_{i,\Gamma'}\in W_i(\Gamma'_i)$, there exist positive constants $C_0$ and $C_1$ independent of $h_i, H_i, \al_i(x)$ and $u_{i,\Gamma'}$, such that
\begin {align}\label {dnormequiv}
C_0 d_i(\H_i u_{i,\Gamma'}, \H_i u_{i,\Gamma'})\leq a'_i(\H'_i u_{i,\Gamma'}, \H'_i u_{i,\Gamma'})\leq C_1 d_i(\H_i u_{i,\Gamma'}, \H_i u_{i,\Gamma'}).
\end {align}
\end {corollary}

Let us introduce the product spaces
\begin {align}\label {Wspace}
W(\om') := \prod_{i=1}^N W_i(\om'_i)\qquad\text{and}\qquad W(\Gamma') := \prod_{i=1}^N W_i(\Gamma'_i).
\end {align}
That is, a function $u\in W(\om')$ means that $u = \{u_i\}_{i=1}^N$ with $u_i\in W_i(\om'_i)$, and a function $u_{\Gamma'} \in W(\Gamma')$ means that $u_{\Gamma'} = \{u_{i,\Gamma'}\}_{i=1}^N$ with $u_{i,\Gamma'}\in W_i(\Gamma'_i)$; see~\eqref {wspacelocal1} and~\eqref {wspacelocal2} for the definitions of $W_i(\om'_i)$ and $W_i(\Gamma'_i)$, and also~\eqref {notation} for notation. We also define
\begin {align}\label {Schurglobal}
S' := \text{diag}\{S'_{1},\cdots,S'_{N}\},
\end {align}
where $S'_i$ is given in~\eqref {Schurlocal}.

%The space $W(\Gamma')$ which was defined only on $\Gamma'$, will also be interpreted below as the subspace of $W(\om')$ of functions which are discrete harmonic in the sense of $\H'_i$ in each $\om_i$, and with the defined boundary values on $\Gamma'_i$.

\subsection {FEIT-DP Problem} Secondly, we formulate~\eqref {disp} as a constrained minimization problem.

With a similar decomposition as~\eqref {wspacelocal1}, we can partition $W_i(\Gamma'_i)$ as \begin {equation}\label {wspaceboundary}
\begin {aligned}
W_i(\Gamma'_i) = W_i(\Gamma_i)\times \prod_{j\in\E_i^0} W_i(\bar{E}_{ji}),
\end {aligned}
\end {equation}
where $W_i(\Gamma_i)$ is the trace of the space $V_h(\om_i)$ on $\Gamma_i$. A function $u_i\in W_i(\Gamma'_i)$ can be written as
\begin {align}\label {IBdecomp}
u_i = \{(u_i)_i,\{(u_i)_j\}_{j\in\E_i^0}\},
\end {align}
where $(u_i)_i$ is $u_i$ restricted to $\bar{E}_{ij}$ and $(u_i)_j$ is $u_i$ restricted to $\bar{E}_{ji}$ for all $j\in \E_i^0$.

We consider $\whw(\Gamma')$ as the subspace of $W(\Gamma')$ which contains the continuous functions on $\Gamma$. A function $u = \{u_i\}_{i=1}^N \in W(\Gamma')$ is defined to be continuous on $\Gamma$ in the sense that for all $1 \leq i\leq N$ we have
\begin {equation}\label {gammacont}
\left\{\begin {aligned}
(u_i)_i(x) &= (u_j)_i(x)\;\;\text{for all}\;\; x\in\bar{E}_{ij}&&\;\;\text{for all}\;\; j\in\E_i^0,\\
(u_i)_j(x) &= (u_j)_j(x)\;\;\text{for all}\;\; x\in\bar{E}_{ji}&&\;\;\text{for all}\;\; j\in\E_i^0.
\end {aligned}\right.
\end {equation}
We say that $u = \{u_i\}_{i=1}^N\in W(\om')$, where $u_i = (u_{i,I}, u_{i,\Gamma})$ with $u_{i,I}\in W_i(I_i)$ and $u_{i,\Gamma}\in W_i(\Gamma'_i)$, is continuous on $\Gamma$ if $\{u_{i,\Gamma}\}_{i=1}^N\in W(\Gamma')$ satisfies the continuity condition~\eqref {gammacont}. The subspace of $W(\om')$ of functions which are continuous on $\Gamma$ is denoted by $\whw(\om')$; c.f., Definition 3.3 in~\cite {DGS13}. Note that there is a one-to-one correspondence between vectors in $V_h(\om)$ and $\whw(\om')$.

%We also denote by $\whw(\Gamma')$ the subspace of $\whw(\om')$ of functions which are discrete harmonic in the sense of $\H'_i$ in each $\om_i$; c.f., Definition 3.3 in~\cite {DGS13}. Note that there is a one-to-one correspondence between vectors in $X(\om)$ and $\whw(\om')$.

Next we define the nodal points associated with the corner variables by
\begin {align}\label {corner}
\V'_i := \V_i \bigcup \{\cup_{j\in\E_i^0} \pa E_{ji}\}\quad\text{where}\quad \V_i:= \{\cup_{j\in\E_i^0} \pa E_{ij}\}.
\end {align}
We now consider the subspace $\wtw(\om') \subset W(\om')$ and $\wtw(\Gamma') \subset W(\Gamma')$ as the space of functions that are continuous on all the $\V'_i$. A function $u = \{u_i\}_{i=1}^N \in W(\Gamma')$ is defined to be continuous at the corners $\V'_i$ in the sense that for all $1 \leq i\leq N$ we have
\begin {equation}\label {cornercont}
\left\{\begin {aligned}
(u_i)_i(x) &= (u_j)_i(x)\;\;\text{at}\;\; x\in\pa E_{ij}&&\;\;\text{for all}\;\; j\in\E_i^0,\\
(u_i)_j(x) &= (u_j)_j(x)\;\;\text{at}\;\; x\in\pa E_{ji}&&\;\;\text{for all}\;\; j\in\E_i^0.
\end {aligned}\right.
\end {equation}
We say that $u = \{u_i\}_{i=1}^N\in W(\om')$, where $u_i = (u_{i,I}, u_{i,\Gamma})$ with $u_{i,I}\in W_i(I_i)$ and $u_{i,\Gamma}\in W_i(\Gamma'_i)$, is continuous on $\V'_i$ if $\{u_{i,\Gamma}\}_{i=1}^N\in W(\Gamma')$ satisfies the continuity condition~\eqref {cornercont}. The subspace of $W(\om')$ of functions which are continuous on $\V'_i$ is denoted by $\wtw(\om')$; c.f., Definition 4.1 in~\cite {DGS13}. Note that $\whw(\Gamma')\subset\wtw(\Gamma')\subset W(\Gamma')$.

%The subspace of $W(\om')$ of continuous functions at the corners $\V'_i$ for all $1 \leq i\leq N$ is denoted by $\wtw(\om')$, and the subspace of $\wtw(\om')$ of functions which are discrete harmonic in the sense of $\H'_i$ is denoted by $\wtw(\Gamma')$; c.f., Definition 4.1 in~\cite {DGS13}. Note that $\whw(\Gamma')\subset\wtw(\Gamma')\subset W(\Gamma')$.

%For any vector $u = \{u_i\}_{i=1}^N \in W(\om')$, we can obtain a vector $u\in\wtw(\om')$ by using the continuity condition~\eqref {cornercont} of $u$ on $\V'_i$ for all $1 \leq i\leq N$.

We can represent $u\in\wtw(\om')$ as $u = (u_I, u_\Pi, u_\Delta)$, where the subscript $I$ refers to the interior degrees of freedom at nodal points $I$; see~\eqref {notation}, the $\Pi$ refers to the primal($\Pi$) variables at the corners $\V'_i$ for all $1 \leq i\leq N$, and the $\Delta$ refers to the dual($\Delta$) variables at the remaining nodal points on $\Gamma'_i\backslash\V'_i$ for all $1 \leq i\leq N$. Similarly, a vector $u\in\wtw(\Gamma')$ can be uniquely decomposed as $u = (u_\Pi, u_\Delta)$. Therefore, we can represent $\wtw(\Gamma') = \whw_\Pi(\Gamma')\times W_\Delta(\Gamma')$, where $\whw_\Pi(\Gamma')$ and $W_\Delta(\Gamma')$ refer to the $\Pi-$ and $\Delta-$degrees of freedom of $\wtw(\Gamma')$, respectively.

Let $\wta$ be the stiffness matrix obtained by restricting the block diagonal matrix $A'$
from $W(\om')$ to $\wtw(\om')$, where $A' = \text{diag}\{A'_1,\cdots,A'_N\}$.
Note that the matrix $\wta$ is no longer block diagonal since there are couplings between primal($\Pi$) variables. Using the decomposition $u = (u_I, u_{II}, u_\Delta)\in\wtw(\om')$, we can partition $\wta$ as
\begin{equation}\label {wta}
\wta=
\left(
\begin{array}{ccc}
A'_{II} & A'_{I\Pi} & A'_{I\Delta} \\
A'_{\Pi I} & \widetilde{A}_{\Pi\Pi} & A'_{\Pi\Delta} \\
A'_{\Delta I} & A'_{\Delta \Pi} & A'_{\Delta \Delta}
\end{array}\right).
\end{equation}
Note that the only coupling across subdomains are through the $\Pi$ variables where the matrix $\wta$ is subassembled.

Once the variables of $I$ and $\Pi$ sets are eliminated, the Schur complement matrix associated with the $\Delta-$variables is obtained of the form
\begin {align}\label {wtschur}
\wts:=A'_{\Delta\Delta}-(A'_{\Delta I}~~ A'_{\Delta\Pi})\left(
                                                            \begin{array}{cc}
                                                              A'_{II} & A'_{I\Pi} \\
                                                              A'_{\Pi I} & \widetilde{A}_{\Pi\Pi} \\
                                                            \end{array}
                                                          \right)^{-1}
\left(
  \begin{array}{c}
    A'_{I\Delta} \\
    A'_{\Pi\Delta} \\
  \end{array}
\right).
\end {align}
Note that $\wts$ is defined on the vector space $W_\Delta(\Gamma')$.

\begin {lemma}\label {lemSchur}
Let $\wta$ and $\wts$ be defined in~\eqref {wta} and~\eqref {wtschur}. For any $u_\Delta\in W_\Delta(\Gamma')$, it holds
\begin {align*}
\langle\wts u_\Delta, u_\Delta\rangle=\min\langle\wta w,w\rangle,
\end {align*}
where the minimum is taken over $w=(w_I, w_\Pi, w_\Delta)\in\wtw(\om')$ with $w_\Delta = u_\Delta$.
\end {lemma}
The proof of the above lemma can be found in Lemma 6.22 of~\cite {TW05} and Lemma 4.2 of~\cite {MT01}.\\

Next we introduce some notations to define the jump matrix $B_\Delta$. The vector space $W_\Delta(\Gamma')$ can be further decomposed as
\begin {align}\label {wdelgamma}
W_\Delta(\Gamma') := \prod_{i=1}^N W_{i,\Delta}(\Gamma'_i),
\end {align}
where the local space $W_{i,\Delta}(\Gamma'_i)$ includes functions associated with variables at the nodal points of $\Gamma'_i\backslash \V'_i$. Hence, a vector $u_\Delta\in W_\Delta(\Gamma')$ can be represented as $u_\Delta = \{u_{i,\Delta}\}_{i=1}^N$ with $u_{i,\Delta}\in W_{i,\Delta}(\Gamma'_i)$. Moreover, the vector $u_{i,\Delta}\in W_{i,\Delta}(\Gamma'_i)$ can be partitioned as
\begin {align*}
u_{i,\Delta} = \{(u_{i,\Delta})_i, \{(u_{i,\Delta})_j\}_{j\in\E_i^0}\}
\end {align*}
with $(u_{i,\Delta})_i=u_{i,\Delta}|_{\Gamma_i\backslash\V_i}$ and
$(u_{i,\Delta})_j=u_{i,\Delta}|_{E_{ji}}$. In order to measure the jump of $u_\Delta\in W_\Delta(\Gamma')$ across the $\Delta-$nodes, we introduce the space
\begin {align*}
\whw_\Delta(\Gamma):=\prod_{i=1}^N V_h(\Gamma_i\backslash\V_i),
\end {align*}
where $V_h(\Gamma_i\backslash\V_i)$ is the restriction of $V_h(\om_i)$ to $\Gamma_i\backslash\V_i$. The jumping matrix $B_\Delta : W_\Delta(\Gamma')\rightarrow\whw_\Delta(\Gamma)$ is constructed as follows: let $u_\Delta = \{u_{i,\Delta}\}_{i=1}^N\in W_\Delta(\Gamma')$ and let $v:= B_\Delta u_\Delta$ where $v = \{v_i\}_{i=1}^N\in\whw_\Delta(\Gamma)$ satisfies
\begin {align}\label {jumpcondition}
v_i = (u_{i,\Delta})_i - (u_{j,\Delta})_i\;\;\text{on}\;\; E_{ijh}\;\;\text{for all}\;\; j\in \E_i^0.
\end {align}
The jumping matrix $B_\Delta$ can be written as
\begin {align}\label {jumpmatrix}
B_\Delta =(B_\Delta^{(1)},B_\Delta^{(2)},\cdots,B_\Delta^{(N)}),
\end {align}
where the rectangular matrix $B_\Delta^{(i)}$ consists of columns of $B_\Delta$ attributed to the $i-$th components of the product space $W_{\Delta}(\Gamma')$. The entries of $B_\Delta^{(i)}$ consist of values of $\{0, 1, -1\}$. It is easy to see that $Range(B_\Delta)=\whw_{\Delta}(\Gamma)$, and $B_\Delta$ has full rank. In addition, if $u = (u_\Pi, u_\Delta)\in\wtw(\Gamma')$ and $B_\Delta u_\Delta = 0$ then $u \in \whw(\Gamma')$.

We can reformulate the discrete problem~\eqref {disp}, on the space of $W_\Delta(\Gamma')$, as a minimization problem with constraints given by the continuity requirement: Find $u^*_\Delta\in W_\Delta(\Gamma')$ such that
\begin {align}\label {miniconst}
\mathcal{J}(u^*_\Delta) = \min\mathcal{J}(v_\Delta),
\end {align}
where the minimum is taken over $v_\Delta\in W_\Delta(\Gamma')$ with constraints $B_\Delta v_\Delta=0$. The objective function
\begin {align}\label {object}
\mathcal{J}(v_\Delta):=\frac{1}{2}\langle\wts v_\Delta, v_\Delta\rangle-\langle\tilde{g}_\Delta, v_\Delta\rangle,
\end {align}
where $\wts$ is defined in~\eqref {wtschur} and
\begin{align*}
\tilde{g}_\Delta:=f_\Delta-(A'_{\Delta I}~~ A'_{\Delta\Pi})\left(
                                                            \begin{array}{cc}
                                                              A'_{II} & A'_{I\Pi} \\
                                                              A'_{\Pi I} & \widetilde{A}_{\Pi\Pi} \\
                                                            \end{array}
                                                          \right)^{-1}
\left(
  \begin{array}{c}
    f_I \\
    f_{\Pi}\\
  \end{array}
\right).
\end{align*}
Here $f=\{f_i\}_{i=1}^N \in V_h(\Omega)$, where $f_i$ is the load vector associated with the subdomain $\om_i$, and $f$ can be represented as $f = (f_I, f_\Pi, f_{\Gamma\backslash\Pi})$. The forcing term $f_\Delta\in W_\Delta(\Gamma')$ is defined by $f_\Delta =\{f_{i,\Delta}\}_{i=1}^N$, where the entries $f_{i,\Delta}$ are defined as $\int_{\om_i} f v_{i,\Delta} dx$ when $v_{i,\Delta}$ are the canonical basis functions of $W_{i,\Delta}(\Gamma'_i)$.

Note that $\wta$ and $\wts$ are both symmetric and positive definite; see also Lemma~\ref {lemSchur}. By introducing a set of Lagrange multipliers $\lambda\in \whw_{\Delta}(\Gamma)$, to enforce the continuity constraints, we obtain the following saddle point formulation of~\eqref {miniconst}: Find $u^*_\Delta\in W_\Delta(\Gamma')$ and $\lambda^*\in \whw_{\Delta}(\Gamma)$ such that
\begin {equation}\label {saddle}
\left\{\begin {aligned}
\wts u^*_\Delta+B_\Delta^T\lambda^*=\tilde{g}_\Delta\\
B_\Delta u^*_\Delta=0.
\end {aligned}\right.
\end {equation}
This reduces to
\begin{align}
F\lambda^*=d,
\end{align}
where
\begin {align}\label {reduce}
F:=B_\Delta\wts^{-1}B_\Delta^T\quad\text{and}\quad d:=B_\Delta\wts^{-1}\tilde{g}_\Delta.
\end {align}
Once $\lambda^*$ is computed, we can back solve and obtain
\begin {align}\label {solution}
u^*_\Delta=\wts^{-1}(\tilde{g}_\Delta-B_\Delta^T\lambda^*).
\end {align}

\subsection {FEIT-DP Preconditioner} We will now define a preconditioner for $F$ in~\eqref {reduce}.

Let us introduce the diagonal scaling matrix $D_\Delta^{(i)}$, which maps $W_{i,\Delta}(\Gamma'_i)$ into itself, for all $1\leq i\leq N$. Each of the diagonal entries of $D_\Delta^{(i)}$ corresponds to one $\Delta-$node, and it is given by the weighted counting function~\cite {PS08}
\begin {align}\label {weightfun}
\delta_j^\dag(x):=\frac{\overline{\al}_j}{\overline{\al}_j + \overline{\al}_i}
\;\;\text{for all}\;\; x\in\{E_{ijh}\cup E_{jih}\}\;\;\text{for all}\;\; j\in\E_i^0,
\end {align}
where $\overline{\al}_i$ is defined in~\eqref {coef}. Note that one edge is shared by two subdomains. The union of all these functions $\delta_j^\dag(x)$ provides a partition of unity on all $\Delta-$nodes.

We also define
\begin {align}\label {BDDelta}
B_{D,\Delta}:= \Lr{B_\Delta^{(1)}D_\Delta^{(1)},\cdots,B_\Delta^{(N)}D_\Delta^{(N)}}.
\end {align}

An important role will be played by the operator
\begin {align}\label {BDproject}
P_\Delta := B_{D,\Delta}^T B_\Delta,
\end {align}
which maps $W_\Delta(\Gamma')$ into itself. It is easy to check that for $w_\Delta = \{w_{i,\Delta}\}_{i=1}^N\in W_\Delta(\Gamma')$ and $v_\Delta:=P_\Delta w_\Delta$, we have
\begin {align}\label {BDprodef1}
(v_{i,\Delta})_i(x)=\delta_j^\dag(x)[(w_{i,\Delta})_i(x)-(w_{j,\Delta})_i(x)]\,\,\text{for all}\;\; x\in E_{ijh},
\end {align}
\begin {align}\label {BDprodef2}
(v_{i,\Delta})_j(x)=\delta_j^\dag(x)[(w_{i,\Delta})_j(x)-(w_{j,\Delta})_j(x)]\,\,\text{for all}\;\; x\in E_{jih},
\end {align}
where $\delta_j^\dag(x)$ is defined in~\eqref {weightfun}. Hence, $P_\Delta$ preserves jumps in the sense that
\begin {align}\label {prejump}
B_\Delta P_\Delta = B_\Delta,
\end {align}
which implies that $P_\Delta$ is a projection with $P_\Delta^2 = P_\Delta$.

Define
\begin {align}\label {SchurDelta}
S'_\Delta := \text{diag}\{S'_{1,\Delta},\cdots,S'_{N,\Delta}\},
\end {align}
where $S'_{i,\Delta}$ is the local Schur complement $S'_i$, see~\eqref {Schurlocal}, restricted to $W_{i,\Delta}(\Gamma'_i)$ from $W_i(\Gamma'_i)$, i.e., $S'_{i,\Delta}$ is obtained from $S'_i$ by deleting rows and columns associated with the variables at nodal points of $\V'_i\subset\Gamma'_i$.

The FETI-DP method is the standard preconditioned conjugate gradient algorithm for solving the preconditioned system
\begin {align*}
M^{-1} F\lambda = M^{-1} d
\end {align*}
with the preconditioner
\begin {align}\label {precondition}
M^{-1}:=B_{D,\Delta}S'_\Delta B_{D,\Delta}^T=\sum_{i=1}^N B_{\Delta}^{(i)}D_{\Delta}^{(i)}S'_{i,\Delta}D_{\Delta}^{(i)}(B_{\Delta}^{(i)})^T.
\end {align}
Note that $M^{-1}$ is a block diagonal matrix and each block is invertible since $S'_{i,\Delta}$ and $D_{\Delta}^{(i)}$ are invertible, and $B_{\Delta}^{(i)}$ has full rank.

\section {Condition Number Estimate for FETI-DP Preconditioner} The main result of our paper is included in the following theorem, which gives an estimate of the condition number for the preconditioned FETI-DP operator $M^{-1}F$.

\begin {theorem}\label {mainthm}
For any $\lambda\in \whw_\Delta(\Gamma)$, there exists a positive constant $C$ independent of $h_i$, $H_i$, $\al(x)$ and $\lambda$ such that
\begin {align}\label {condition}
\langle M\lambda, \lambda\rangle\leq \langle F\lambda, \lambda\rangle\leq C \beta\langle M\lambda, \lambda\rangle,
\end {align}
where
\begin {align}\label {beta1}
\beta = (\frac{H}{h})^2 \max_{i=1}^N \frac{\overline{\al}_i}{\underline{\al}_i}
\end {align}
with $H/h=\max_{i=1}^N H_i/h_i$. If for any $1\leq i \leq N$ the coefficient $\al(x)$ in the subdomain $\om_i$ satisfies
\begin {align}\label {coefcond}
\al(x)\geq\underline{\al}_i\qquad\text{for all}\qquad x\in\om_i,
\end {align}
then we have
\begin {align}\label {beta2}
\beta = \frac{H}{h}(1+\log \frac{H}{h})^2 \max_{i=1}^N  \frac{\overline{\al}_i}{\underline{\al}_i} .
\end {align}
\end {theorem}
\begin {proof} By the general abstract theory for FETI-DP method, see~\cite {MT01} and Theorem 6.35 of~\cite {TW05}, the proof of the lower and upper bound in~\eqref {condition} follows by checking Lemma~\ref {lowerbound} and Lemma~\ref {upperbound} as below, respectively.
\end {proof}

For clarity, we will use the following norms for $w=(w_\Pi, w_\Delta)\in\wtw(\Gamma')$ with $w_\Delta\in W_\Delta(\Gamma')$:
\begin {align*}
\|w\|^2_{S'} := \langle S'w, w\rangle,\quad\quad \|w_\Delta\|^2_{\wts} := \langle \wts w_\Delta, w_\Delta\rangle,
\end {align*}
and
\begin {align}\label {normsp}
\|w_\Delta\|^2_{S'_\Delta} := \langle S'_\Delta w_\Delta, w_\Delta\rangle = \langle S' \left(
                                                            \begin{array}{c}
                                                              0 \\
                                                              w_\Delta
                                                            \end{array}
                                                          \right),
                                                          \left(
                                                            \begin{array}{c}
                                                              0 \\
                                                              w_\Delta
                                                            \end{array}
                                                          \right) \rangle,
\end {align}
where $S'$, $\wts$ and $S'_\Delta$ are defined in~\eqref {Schurglobal},~\eqref {wtschur} and~\eqref {SchurDelta}, respectively.

\begin {lemma}\label {lowerbound}
For any $\mu\in \whw_\Delta(\Gamma)$ there exists a $w_\Delta\in W_\Delta(\Gamma')$ such that
\begin {align*}
\mu = B_\Delta w_\Delta
\end {align*}
with
\begin {align*}
P_\Delta w_\Delta = w_\Delta
\end {align*}
and
\begin {align*}
\|w_\Delta\|_{\wts} \leq \|P_\Delta w_\Delta\|_{S'_\Delta}.
\end {align*}
\end {lemma}
\begin {proof} For any $\mu\in \whw_\Delta(\Gamma)$, there exists an element $v_\Delta\in W_\Delta(\Gamma')$ such that
\[
\mu = B_\Delta v_\Delta,
\]
since $B_\Delta$ has full rank.

Note that $P_\Delta$ is a projection which maps $W_\Delta(\Gamma')$ to itself. By choosing
\[
w_\Delta = P_\Delta v_\Delta\in W_\Delta(\Gamma'),
\]
we can easily obtain
\begin {align*}
P_\Delta w_\Delta = P^2_\Delta v_\Delta= P_\Delta v_\Delta=w_\Delta,
\end {align*}
and
\begin {align*}
B_\Delta w_\Delta = B_\Delta P_\Delta v_\Delta = B_\Delta v_\Delta=\mu,
\end {align*}
where we have used~\eqref {prejump}.

It follows from Lemma~\ref {lemSchur} that
\begin {align*}
\|w_\Delta\|^2_{\wts} = \min \langle\wta v, v\rangle \leq \min\langle \wta\hat{v}, \hat{v}\rangle = \|w_\Delta\|^2_{S'_\Delta} = \|P_\Delta w_\Delta\|^2_{S'_\Delta},
\end {align*}
where the first minimum is taken over $v = (v_I, v_\Pi, v_\Delta)\in\wtw(\om')$ with $v_\Delta = w_\Delta$, and the second one over $\hat{v} = (\hat{v}_I, \hat{v}_\Pi, \hat{v}_\Delta)$ with $\hat{v}_\Pi = 0$ and $\hat{v}_\Delta = w_\Delta$.
\end {proof}

The next lemma gives us a crucial estimate of the norm of $P_\Delta$.

\begin {lemma}\label {upperbound}
For any $w_\Delta\in W_\Delta(\Gamma')$ it holds that
\begin {align}\label {normofP}
\|P_\Delta w_\Delta\|^2_{S'_\Delta}\leq C\beta \|w_\Delta\|^2_{\wts},
\end {align}
where $\beta$ is defined in~\eqref {beta1} or/and~\eqref {beta2}, and $C$ is a positive constant independent of $h_i$, $H_i$, $\alpha(x)$ and $w_\Delta$.
\end {lemma}
\begin {proof}
For any $w_\Delta\in W_\Delta(\Gamma')$, let $w = (w_\Pi, w_\Delta)\in\wtw(\Gamma')$ be the solution of
\begin {align}\label {righthand}
\langle S' w, w\rangle = \min \langle S' v, v\rangle = \langle \wts w_\Delta, w_\Delta\rangle,
\end {align}
where the minimum is taken over $v = (v_\Pi, v_\Delta)\in \wtw(\Gamma')$ with $v_\Pi\in\whw_\Pi(\Gamma')$ and $v_\Delta = w_\Delta$.

%Therefore, the term $\|w_\Delta\|_{\wts}$ in the right hand side of~\eqref {normofP} can be replaced by $\|w\|_{S'}$.

We can represent $w$ as $w = \{w_i\}_{i=1}^N \in\wtw(\Gamma')$ with $w_i\in W_i(\Gamma'_i)$. We define linear functions to approximate $w_i$ on $\bar{E}_{ij}$ and $\bar{E}_{ji}$ with $j\in\E_i^0$ as follows:
\begin {align*}
I_{E_{ij}} (w_i)_i(x) \;\;\text{is linear on}\;\;\bar{E}_{ij}\;\;\text{and}\;\;I_{E_{ij}} (w_i)_i(x) = (w_i)_i(x)\;\;\text{for all}\;\; x\in\pa E_{ij},
\end {align*}
and
\begin {align*}
I_{E_{ji}} (w_i)_j(x) \;\;\text{is linear on}\;\;\bar{E}_{ji}\;\;\text{and}\;\;I_{E_{ji}} (w_i)_j(x) = (w_i)_j(x)\;\;\text{for all}\;\; x\in\pa E_{ji}.
\end {align*}
Let $\hat{w}=\{\hat{w}_i\}_{i=1}^N$ with $\hat{w}_i\in W_i(\Gamma'_i)$ be defined by
\begin {align*}
(\hat{w}_i)_i(x) = I_{E_{ij}} (w_i)_i(x) \;\;\text{for all}\;\; x\in\bar{E}_{ijh}\;\;\text{for all}\;\; j\in\E_i^0,
\end {align*}
and
\begin {align*}
(\hat{w}_i)_j(x) = I_{E_{ji}} (w_i)_j(x) \;\;\text{for all}\;\; x\in\bar{E}_{jih}\;\;\text{for all}\;\; j\in\E_i^0.
\end {align*}
Note that $\hat{w}\in\whw(\Gamma')$; see~\eqref {gammacont}. Therefore, representing $\hat{w} = (\hat{w}_\Pi, \hat{w}_\Delta)$, we have $B_\Delta \hat{w}_\Delta = 0$. Using the definition of $P_\Delta$~\eqref {BDproject}, we have
\begin {align*}
P_\Delta w_\Delta = B^T_{D,\Delta} B_\Delta w_\Delta = B^T_{D,\Delta} B_\Delta (w_\Delta - \hat{w}_\Delta) =  P_\Delta (w_\Delta - \hat{w}_\Delta).
\end {align*}

Define $v\in\wtw(\Gamma')$ to be equal to $P_\Delta(w_\Delta - \hat{w}_\Delta)$ at the $\Delta-$nodes, and equal to zero at the $\Pi-$nodes. Let us represent $v = \{v_i\}_{i=1}^N$ with $v_i\in W_i(\Gamma'_i)$ and
\begin {align*}
v_i = \{(v_i)_i,\{(v_i)_j\}_{j\in\E_i^0}\},
\end {align*}
where $(v_i)_i\in W_i(\Gamma_i)$; see~\eqref {wspaceboundary} and~\eqref {IBdecomp}. Using~\eqref {BDprodef1} and~\eqref {BDprodef2}, it is easy to check that
\begin {align}\label {vBDprodef1}
(v_i)_i = \delta_j^\dag(x) [(w_i-\hat{w}_i)_i - (w_j-\hat{w}_j)_i],
\end {align}
and
\begin {align}\label {vBDprodef2}
(v_i)_j = \delta_j^\dag(x) [(w_i-\hat{w}_i)_j - (w_j-\hat{w}_j)_j].
\end {align}
We denote by $V_h(\pa\om_i)$ the space of continuous and piecewise linear functions on the local boundaries $\pa\om_i$. It is obvious that $(v_i)_i\in V_h(\pa\om_i)$.
By the definitions of $S'_\Delta$ and $S'$, ~\eqref {SchurDelta},~\eqref {Schurglobal}, and~\eqref {normsp}, we have
\begin {align*}
\|P_\Delta w_\Delta\|^2_{S'_\Delta} = \|v\|^2_{S'} =\sum_{i=1}^N \|v_i\|^2_{S'_i},
\end {align*}
where
\begin {align*}
\|v_i\|^2_{S'_i} = \langle S'_i v_i, v_i\rangle = a'_i(\H'_i v_i, \H'_i v_i)
\end {align*}
with the discrete harmonic extension $\H'_i$ defined in~\eqref {dishar}.

With~\eqref {righthand}, to prove~\eqref {normofP}, we need to show that
\begin {align*}
\sum_{i=1}^N a'_i(\H'_i v_i, \H'_i v_i) \leq C\beta \sum_{i=1}^N a'_i(\H'_i w_i, \H'_i w_i).
\end {align*}
By Corollary~\ref {coro1} it remains to prove
\begin {align*}
 \sum_{i=1}^N d_i(\H_i v_i, \H_i v_i) \leq C\beta \sum_{i=1}^N d_i(\H_i w_i, \H_i w_i),
\end {align*}
with
\begin {equation}\label {estimate}
\begin {aligned}
%d_i(\H_i v_i, \H_i v_i)&=\sum_{\tau\subset\om_i}\alpha_\tau||\na (\H_i v_i)||_{L^2(\tau)}^2+\sum_{j\in\E_i}\sum_{e\subset E_{ij}}\frac{\delta\al_e}{l_{ij}h_{ij}}||(v_i)_i-(v_i)_j||_{L^2(e)}^2,\\
%&= I_1+I_2.
d_i(\H_i v_i, \H_i v_i)&=||\al_i^{1/2}\na (\H_i v_i)||_{L^2(\om_i)}^2+\sum_{j\in\E_i}\frac{1}{l_{ij}}\frac{\delta}{h_{ij}}||\al_{ij}^{1/2}[(v_i)_i-(v_i)_j]||_{L^2(E_{ij})}^2\\
&= I_1+I_2.
\end {aligned}
\end {equation}

First we consider the term $I_2$ of~\eqref {estimate}. For $j\in\E_i^\pa$ the proof is trivial due to the specific choices of parameters. For $j\in\E_i^0$, it follows from~\eqref {vBDprodef1} and~\eqref {vBDprodef2} that
\begin {align*}
||(v_i)_i-(v_i)_j||_{L^2(e)}^2&= (\delta_j^\dag(x))^2\|(w_i-\hat{w}_i)_i-(w_j-\hat{w}_j)_i-(w_i-\hat{w}_i)_j+(w_j-\hat{w}_j)_j\|^2_{L^2(e)}\\
&\leq \|(w_i)_i-(w_i)_j\|^2_{L^2(e)}+\|(w_j)_i-(w_j)_j\|^2_{L^2(e)},
\end {align*}
since $\hat{w}\in\whw(\Gamma')$, and $\delta_j^\dag\in (0, 1)$. Here $e$ is a fine edge on the mesh that is obtained by merging $\T_h(\om_i)$ and $\T_h(\om_j)$ along $E_{ij}$. We recall that $\al_{ij}$ is constant on each $e\subset E_{ij}$ and denoted by $\al_{ij}^e$.

By summing up, we finally get
\begin {equation}
\begin {aligned}\label {estimateofI2}
I_2 &\leq C \sum_{j\in\E_i}\frac{1}{l_{ij}}\frac{\delta}{h_{ij}}\sum_{e\subset E_{ij}}\al_{ij}^e\|(w_i)_i-(w_i)_j\|^2_{L^2(e)}+\\
&\quad +C \sum_{j\in\E_i}\frac{1}{l_{ij}}\frac{\delta}{h_{ij}}\sum_{e\subset E_{ij}}\al_{ij}^e\|(w_j)_i-(w_j)_j\|^2_{L^2(e)}\\
&=C\sum_{j\in\E_i}\frac{1}{l_{ij}}\frac{\delta}{h_{ij}}\|\al_{ij}^{1/2}[(w_i)_i-(w_i)_j]\|^2_{L^2(E_{ij})}+\\
&\quad+C \sum_{j\in\E_i}\frac{1}{l_{ij}}\frac{\delta}{h_{ij}}\|\al_{ij}^{1/2}[(w_j)_i-(w_j)_j]\|^2_{L^2(E_{ij})}\\
& \leq C \{d_i(\H_i w_i, \H_i w_i) + \sum_{j\in\E_i} d_j(\H_j w_j, \H_j w_j)\}.
\end {aligned}
\end {equation}

Now we estimate the first term $I_1$ of~\eqref {estimate}. Here we introduce two semi-norms defined on $\Gamma_i$ as follows: for any $u_i \in W_i(\Gamma_i)$,
\begin {align}\label {semiinter1}
%|u_i|_{B_i}^2:=\min \{\sum_{\tau\subset \Omega_i}\al_\tau\|\na\tilde{u}_i\|^2_{L^2(\tau)}\;\;:\;\; \tilde{u}_i \in W_i(\om_i)\;\;\text{and}\;\; \tilde{u}_i=u\;\;\text{on}\;\; \pa\om_i\},
|u_i|_{B_i}^2:=\min \{\|\al_i^{1/2}\na\tilde{u}_i\|^2_{L^2(\om_i)}\;\;:\;\; \tilde{u}_i \in V_h(\om_i)\;\;\text{and}\;\; \tilde{u}_i|_{\pa\om_i}=u_i\},
\end {align}
and
\begin {align}\label {semiinter2}
%|u_i|_{H^{1/2}(\pa\om_i)}^2:=\min \{\|\na\tilde{u}_i\|^2_{L^2(\om_i)}\;\;:\;\; \tilde{u}_i \in W_i(\om_i)\;\;\text{and}\;\; \tilde{u}_i=u\;\;\text{on}\;\; \pa\om_i\}.
|u_i|_{H^{1/2}(\pa\om_i)}^2:=\min \{\|\na\tilde{u}_i\|^2_{L^2(\om_i)}\;\;:\;\; \tilde{u}_i \in V_h(\om_i)\;\;\text{and}\;\; \tilde{u}_i|_{\pa\om_i}=u_i\}.
\end {align}
%Define the discrete harmonic extension $\accentset{\circ}{\H}_i: W_i(\Gamma_i) \rightarrow W_i(\om_i)$, which satisfies
%\begin {align}
%a_i(\accentset{\circ}{\H}_i u_i, \tilde{u}_i) = 0\quad\text{for all}\quad \tilde{u}_i\in\accentset{\circ}{W}_i(\om_i)
%\end {align}
%with $\accentset{\circ}{\H}_i u_i = u_i$ on $\Gamma_i$.
Denote by $\accentset{\circ}{\H}_i: W_i(\Gamma_i) \rightarrow V_h(\om_i)$ as the discrete harmonic extension in the sense of $a_i(\cdot,\cdot)$. Hence the function $\accentset{\circ}{\H}_i u_i$ is the minimizing function of~\eqref {semiinter1}.

Note that $\H_i v_i$ at the interior nodes depends only on the nodal values of $v_i$ on $\Gamma_i$, i.e., $\accentset{\circ}{\H}_i (v_i)_i = \H_i v_i$ in the interior of subdomains $\om_i$. This implies that
%\begin {equation}\label {I1step1}
%\begin {aligned}
%%I_1 &= \sum_{\tau\subset\om_i}\alpha_\tau||\na [\accentset{\circ}{\H}_i (v_i)_i]||_{L^2(\tau)}^2\\
%%&= |(v_i)_i|_{B_i}^2\\
%%&\leq C \overline{\al}_i \Lr{|(v_i)_i|_{H^{1/2}(\pa\om_i)}^2 + \frac{1}{\eta_i}|(v_i)_i|_{L^2(\pa\om_i)}^2},
%I_1 &= ||\al_i^{1/2}\na [\accentset{\circ}{\H}_i (v_i)_i]||_{L^2(\om_i)}^2\\
%&= |(v_i)_i|_{B_i}^2\\
%&\leq C \overline{\al}_i \Lr{|(v_i)_i|_{H^{1/2}(\pa\om_i)}^2 + \frac{1}{\eta_i}\|(v_i)_i\|_{L^2(\pa\om_i)}^2},
%\end {aligned}
%\end {equation}
\begin {equation}\label {I1step1}
\begin {aligned}
%I_1 &= \sum_{\tau\subset\om_i}\alpha_\tau||\na [\accentset{\circ}{\H}_i (v_i)_i]||_{L^2(\tau)}^2\\
%&= |(v_i)_i|_{B_i}^2\\
%&\leq C \overline{\al}_i \Lr{|(v_i)_i|_{H^{1/2}(\pa\om_i)}^2 + \frac{1}{\eta_i}|(v_i)_i|_{L^2(\pa\om_i)}^2},
I_1 &= ||\al_i^{1/2}\na [\accentset{\circ}{\H}_i (v_i)_i]||_{L^2(\om_i)}^2\\
&= |(v_i)_i|_{B_i}^2\\
&\leq C \overline{\al}_i \Lr{|(v_i)_i|_{H^{1/2}(\pa\om_i)}^2 + \frac{1}{h_i}\|(v_i)_i\|_{L^2(\pa\om_i)}^2},
\end {aligned}
\end {equation}
where we have used the second inequality of Lemma 4.1 in~\cite {PS08}.

We can write $(v_i)_i$ as
\begin {equation}\label {parofuni}
\begin {aligned}
(v_i)_i &= \sum_{j\in\E_i} I^h (\theta_{E_{ij}} (v_i)_i)\\
&=\sum_{j\in\E_i} I^h \Lr{\theta_{E_{ij}} \delta_j^\dag(x) [(w_i-\hat{w}_i)_i - (w_j-\hat{w}_j)_i]},
\end {aligned}
\end {equation}
where $I^h$ is the usual Lagrange interpolation operator, and for $j\in\E_i$ the finite element cut-off function $\theta_{E_{ij}}(x)$ equals to 1 for all $x\in E_{ijh}$ and vanishes on all the other nodes; see Definition 4.2 of~\cite {PS08}.

By~\eqref {weightfun}, we know that
\begin {align*}
\overline{\al}_i\Lr{\delta_j^\dag(x)}^2\leq C \min(\overline{\al}_i, \overline{\al}_j).
\end {align*}
Putting~\eqref {parofuni} into~\eqref {I1step1}, we obtain
%\begin {equation}\label {app2I1step2}
%\begin {aligned}
%I_1 &\leq C \min(\overline{\al}_i, \overline{\al}_j) \sum_{j\in\E_i}\{|\psi_{ij}|_{H^{1/2}(\pa\om_i)}^2 + \frac{1}{\eta_i}|\psi_{ij}|_{L^2(\pa\om_i)}^2 \}\\
%&\leq C \min(\overline{\al}_i, \overline{\al}_j) \sum_{j\in\E_i}\{|a_{ij}|_{H^{1/2}(\pa\om_i)}^2 + |b_{ij}|_{H^{1/2}(\pa\om_i)}^2+\\ &\quad+\frac{1}{\eta_i}\|a_{ij}\|_{L^2(\pa\om_i)}^2  + \frac{1}{\eta_i}\|b_{ij}\|_{L^2(\pa\om_i)}^2\},
%\end {aligned}
%\end {equation}
\begin {equation}\label {app2I1step2}
\begin {aligned}
I_1 &\leq C \min(\overline{\al}_i, \overline{\al}_j) \sum_{j\in\E_i}\{|\psi_{ij}|_{H^{1/2}(\pa\om_i)}^2 + \frac{1}{h_i}|\psi_{ij}|_{L^2(\pa\om_i)}^2 \}\\
&\leq C \min(\overline{\al}_i, \overline{\al}_j) \sum_{j\in\E_i}\{|a_{ij}|_{H^{1/2}(\pa\om_i)}^2 + |b_{ij}|_{H^{1/2}(\pa\om_i)}^2+\\ &\quad+\frac{1}{h_i}\|a_{ij}\|_{L^2(\pa\om_i)}^2  + \frac{1}{h_i}\|b_{ij}\|_{L^2(\pa\om_i)}^2\},
\end {aligned}
\end {equation}
where
\begin {align*}
\psi_{ij} &= I^h \Lr{\theta_{E_{ij}}[(w_i-\hat{w}_i)_i - (w_j-\hat{w}_j)_i]}\\
&= I^h \Lr{\theta_{E_{ij}}[(w_i-\hat{w}_i)_i - (w_j-\hat{w}_j)_j]} + I^h \Lr{\theta_{E_{ij}}[(w_j-\hat{w}_j)_j - (w_j-\hat{w}_j)_i]}\\
&:= a_{ij} + b_{ij}.
\end {align*}
As stated in~\cite {PS08} that $|I^h (\theta_{E_{ij}}(w_j-\hat{w}_j)_j)|_{H^{1/2}(\pa\om_i)} \simeq |I^h (\theta_{E_{ij}}(w_j-\hat{w}_j)_j)|_{H^{1/2}(\pa\om_j)}$, since the discrete harmonic extensions from $E_{ij}$ to $\om_i$ and $\om_j$ are equivalent in the corresponding $H^1-$seminorms. Here we refer to Lemma 4.19 of~\cite {TW05} with the two dimensional case, and have
\begin {equation}\label {app2esta21}
\begin {aligned}
|a_{ij}|_{H^{1/2}(\pa\om_i)}^2 &= |I^h \Lr{\theta_{E_{ij}}[(w_i-\hat{w}_i)_i - (w_j-\hat{w}_j)_j]}|_{H^{1/2}(\pa\om_i)}^2\\
&\leq C\Lr{|I^h\Lr{\theta_{E_{ij}}(w_i-\hat{w}_i)_i}|^2_{H^{1/2}(\pa\om_i)}+|I^h\Lr{\theta_{E_{ij}}(w_j-\hat{w}_j)_j}|^2_{H^{1/2}(\pa\om_j)}}\\
&\leq C\Lr{\frac{1}{h_i}\|(w_i-\hat{w}_i)_i\|^2_{L^2(E_{ij})}+\frac{1}{h_j}\|(w_j-\hat{w}_j)_j\|^2_{L^2(E_{ij})}}.
\end {aligned}
\end {equation}
Since $\hat{w}_i$ is a convex combination of the values of $w_i$ at the end points of $E_{ij}$, we can employ the generalized discrete Sobolev inequality, c.f. Lemma 4.5 of~\cite {PS08}, and obtain
%\begin {equation}\label {app2esta22}
%\begin {aligned}
%&\quad\min(\overline{\al}_i, \overline{\al}_j) \sum_{j\in\E_i}|a_{ij}|_{H^{1/2}(\pa\om_i)}^2\\
%&\leq C \sum_{j\in\E_i}\sum_{k=i,j}\overline{\al}_k  \frac{H_k^2}{h_k\eta_k}(1+\log\frac{\eta_k}{h_k}) |\accentset{\circ}{\H}_k (w_k)_k|_{H^1(\om_{k,\eta_k})}^2\\
%&\leq C\sum_{j\in\E_i}\sum_{k=i,j}\frac{\overline{\al}_k}{\underline{\al}_k}\frac{H_k^2}{h_k\eta_k}(1+\log \frac{\eta_k}{h_k}) |(w_k)_k|_{B_k}^2\\
%&\leq C \frac{H^2}{h\eta}(1+\log\frac{\eta}{h})\sum_{j\in\E_i}\sum_{k=i,j}\frac{\overline{\al}_k}{\underline{\al}_k}d_k(\H_k w_k, \H_k w_k).
%\end {aligned}
%\end {equation}
\begin {equation}\label {app2esta22}
\begin {aligned}
&\quad\min(\overline{\al}_i, \overline{\al}_j) \sum_{j\in\E_i}|a_{ij}|_{H^{1/2}(\pa\om_i)}^2\\
&\leq C \sum_{j\in\E_i}\sum_{k=i,j}\overline{\al}_k  \Lr{\frac{H_k}{h_k}}^2 |\accentset{\circ}{\H}_k (w_k)_k|_{H^1(\om_k^h)}^2\\
&\leq C\sum_{j\in\E_i}\sum_{k=i,j}\frac{\overline{\al}_k}{\underline{\al}_k}\Lr{\frac{H_k}{h_k}}^2|(w_k)_k|_{B_k}^2\\
&\leq C \Lr{\frac{H}{h}}^2\sum_{j\in\E_i}\sum_{k=i,j}\frac{\overline{\al}_k}{\underline{\al}_k}d_k(\H_k w_k, \H_k w_k).
\end {aligned}
\end {equation}
With the same argument as~\eqref {app2esta21}, we get
\begin {equation}\label {app2estb1}
\begin {aligned}
\sum_{j\in\E_i}|b_{ij}|_{H^{1/2}(\pa\om_i)}^2
\leq C\sum_{j\in\E_i}\frac{1}{h_j}\|(w_j-\hat{w}_j)_j - (w_j-\hat{w}_j)_i\|^2_{L^2(E_{ij})}.
\end {aligned}
\end {equation}
Let $Q_i$ be the $L_2$ projection on $V_h(E_{ij})$, the restriction of $V_h(\partial \om_i)$ to $\bar{E}_{ij}$ with $h_i-$triangulation on $E_{ij}$. Using the inverse inequality, and the $L_2$ stability of the $L_2$ projection we have
%\begin {equation}\label {app2esta41}
%\begin {aligned}
%&\quad\|(w_j-\hat{w}_j)_j - (w_j-\hat{w}_j)_i\|^2_{L^2(E_{ij})}\\
%&\leq C\{\|Q_i[(w_j)_j - (w_j)_i]\|^2_{L^2(E_{ij})} + \|Q_i(w_j-\hat{w}_j)_j\|^2_{L^2(E_{ij})}+\\
%&\quad+\|(w_j-\hat{w}_j)_j\|^2_{L^2(E_{ij})}+\|(\hat{w}_j)_i - (\hat{w}_j)_j\|^2_{L^2(E_{ij})}\}\\
%&\leq C\{\|(w_j)_j - (w_j)_i\|^2_{L^2(E_{ij})}+\|(w_j-\hat{w}_j)_j\|^2_{L^2(E_{ij})}+\|(\hat{w}_j)_i - (\hat{w}_j)_j\|^2_{L^2(E_{ij})}\}\\
%&\leq C\{\|(w_j)_j - (w_j)_i\|^2_{L^2(E_{ij})} + \frac{H_j^2}{\eta_j}(1+\log\frac{\eta_j}{h_j})|\accentset{\circ}{\H}_j (w_j)_j|_{H^1(\om_{j,\eta_j})}^2+\\
%&+ H_i \max_{\partial E_{ij}}((w_j)_i - (w_j)_j)^2\},
%\end {aligned}
%\end {equation}
\begin {equation}\label {app2esta41}
\begin {aligned}
&\quad\|(w_j-\hat{w}_j)_j - (w_j-\hat{w}_j)_i\|^2_{L^2(E_{ij})}\\
&\leq C\{\|Q_i[(w_j)_j - (w_j)_i]\|^2_{L^2(E_{ij})} + \|Q_i(w_j-\hat{w}_j)_j\|^2_{L^2(E_{ij})}+\\
&\quad+\|(w_j-\hat{w}_j)_j\|^2_{L^2(E_{ij})}+\|(\hat{w}_j)_i - (\hat{w}_j)_j\|^2_{L^2(E_{ij})}\}\\
&\leq C\{\|(w_j)_j - (w_j)_i\|^2_{L^2(E_{ij})}+\|(w_j-\hat{w}_j)_j\|^2_{L^2(E_{ij})}+\|(\hat{w}_j)_i - (\hat{w}_j)_j\|^2_{L^2(E_{ij})}\}\\
&\leq C\{\|(w_j)_j - (w_j)_i\|^2_{L^2(E_{ij})} + \frac{H_j^2}{h_j}|\accentset{\circ}{\H}_j (w_j)_j|_{H^1(\om_j^h)}^2+ H_j \max_{\partial E_{ij}}((w_j)_i - (w_j)_j)^2\},
\end {aligned}
\end {equation}
since $(\hat{w}_j)_i$ and $(\hat{w}_j)_j$ are linear on $E_{ij}$ and $E_{ji}$, respectively. Let $(\bar{w}_j)_j$ be the average of $(w_j)_j$ on $E_{ji}$. By (4.42) in~\cite {DGS13} we obtain
%\begin {equation}\label {app2esta42}
%\begin {aligned}
%&\quad\max_{\pa E_{ij}}((w_j)_i - (w_j)_j)^2\\
%&\leq C\{\frac{1}{h_i}\|(w_j)_i - (w_j)_j\|^2_{L^2(E_{ji})} + \max_{\pa E_{ij}}(Q_i(w_j-\bar{w}_j)_j)^2 + \max_{\pa E_{ji}}((w_j-\bar{w}_j)_j)^2\}\\
%&\leq C\{\frac{1}{h_i}\|(w_j)_i - (w_j)_j\|^2_{L^2(E_{ji})}+ \frac{H_j}{\eta_j}(1+\log\frac{\eta_j}{h_j})|\accentset{\circ}{\H}_j (w_j)_j|_{H^1(\om_{j,\eta_j})}^2\},
%\end {aligned}
%\end {equation}
\begin {equation}\label {app2esta42}
\begin {aligned}
&\quad\max_{\pa E_{ij}}((w_j)_i - (w_j)_j)^2\\
&\leq C\{\frac{1}{h_i}\|(w_j)_i - (w_j)_j\|^2_{L^2(E_{ji})} + \max_{\pa E_{ij}}(Q_i(w_j-\bar{w}_j)_j)^2 + \max_{\pa E_{ji}}((w_j-\bar{w}_j)_j)^2\}\\
&\leq C\{\frac{1}{h_i}\|(w_j)_i - (w_j)_j\|^2_{L^2(E_{ji})}+ \frac{H_j}{h_j}|\accentset{\circ}{\H}_j (w_j)_j|_{H^1(\om_j^h)}^2\},
\end {aligned}
\end {equation}
where we have used (4.10) in~\cite {PS08}, and the $H^{1/2}$ stability of the $L^2$ projection. Substituting~\eqref {app2esta42} into~\eqref {app2esta41}, we have
\begin {equation}\label {app2estatemp}
\begin {aligned}
&\quad\|(w_j-\hat{w}_j)_j - (w_j-\hat{w}_j)_i\|^2_{L^2(E_{ij})}\\
&\leq C\{H_j\frac{1}{h_{ij}}\|(w_j)_i - (w_j)_j\|^2_{L^2(E_{ji})}+\frac{H_j^2}{h_j}|\accentset{\circ}{\H}_j (w_j)_j|_{H^1(\om_j^h)}^2\},
\end {aligned}
\end {equation}
where we have used~\eqref {haraver}. Putting the above inequality into~\eqref {app2estb1}, we obtain
\begin {equation}\label {app2esta43}
\begin {aligned}
&\quad\min(\overline{\al}_i, \overline{\al}_j) \sum_{j\in\E_i}|b_{ij}|_{H^{1/2}(\pa\om_i)}^2\\
&\leq C\sum_{j\in\E_i}\{\frac{H_j}{h_j}\max(\frac{\overline{\al}_i} {\underline{\al}_i}, \frac{\overline{\al}_j} {\underline{\al}_j})\sum_{e\subset E_{ij}}\frac{\al_{ij}^e}{h_{ij}}\|(w_j)_i - (w_j)_j\|^2_{L^2(e)}+\Lr{\frac{H_j}{h_j}}^2\frac{\overline{\al}_j} {\underline{\al}_j}|(w_j)_j|_{B_j}^2\}\\
&\leq C\sum_{j\in\E_i}\{\frac{H_j}{h_j}\max(\frac{\overline{\al}_i} {\underline{\al}_i}, \frac{\overline{\al}_j} {\underline{\al}_j})\frac{1}{h_{ij}}\|\al_{ij}^{1/2}(w_j)_i - (w_j)_j\|^2_{L^2(E_{ij})}+\Lr{\frac{H_j}{h_j}}^2\frac{\overline{\al}_j} {\underline{\al}_j}|(w_j)_j|_{B_j}^2\}\\
&\leq C \Lr{\frac{H}{h}}^2\max_{i=1}^N \frac{\overline{\al}_i} {\underline{\al}_i} \sum_{j\in\E_i}d_j(\H_j w_j, \H_j w_j),
\end {aligned}
\end {equation}
where we used the fact that for all $e\subset E_{ij}$
\[
\min(\overline{\al}_i, \overline{\al}_j) \leq \max(\frac{\overline{\al}_i} {\underline{\al}_i}, \frac{\overline{\al}_j} {\underline{\al}_j}) \min(\underline{\al}_i, \underline{\al}_j) \leq \max(\frac{\overline{\al}_i} {\underline{\al}_i}, \frac{\overline{\al}_j} {\underline{\al}_j}) \al_{ij}^e.
\]
Using the $L^2$ continuity of the nodal interpolation operator $I^h$, and proceeding the same lines of~\eqref {app2esta22}, we have
%\begin {equation}\label {app2esta44}
%\begin {aligned}
%\|a_{ij}\|_{L^2(\pa\om_i)}^2
%&\leq \|(w_i - \hat{w}_i)_i\|^2_{L^2(E_{ij})}+\|(w_j - \hat{w}_j)_j\|^2_{L^2(E_{ij})}\\
%&\leq C \sum_{k=i,j} \frac{H_k^2}{\eta_k}(1+\log\frac{\eta_k}{h_k})|\accentset{\circ}{\H}_k (w_k)_k|_{H^1(\om_{k,\eta_k})}^2.
%\end {aligned}
%\end {equation}
%Hence,
%\begin {equation}\label {app2esta45}
%\begin {aligned}
%\min(\overline{\al}_i, \overline{\al}_j) \sum_{j\in\E_i}\frac{1}{\eta_i}\|a_{ij}\|_{L^2(\pa\om_i)}^2
%\leq C\frac{H^2}{\eta^2}(1+\log\frac{\eta}{h})\sum_{j\in\E_i}\sum_{k=i,j}\frac{\overline{\al}_k} {\underline{\al}_k}d_k(\H_k w_k, \H_k w_k).
%\end {aligned}
%\end {equation}
\begin {equation}\label {app2esta45}
\begin {aligned}
&\quad\min(\overline{\al}_i, \overline{\al}_j) \sum_{j\in\E_i}\frac{1}{h_i}\|a_{ij}\|_{L^2(\pa\om_i)}^2\\
&\leq C\min(\overline{\al}_i, \overline{\al}_j)\sum_{j\in\E_i}\frac{1}{h_i}(\|(w_i - \hat{w}_i)_i\|^2_{L^2(E_{ij})}+\|(w_j - \hat{w}_j)_j\|^2_{L^2(E_{ij})})\\
&\leq C\Lr{\frac{H}{h}}^2\sum_{j\in\E_i}\sum_{k=i,j}\frac{\overline{\al}_k} {\underline{\al}_k}d_k(\H_k w_k, \H_k w_k),
\end {aligned}
\end {equation}
and
\begin {equation}\label {app2esta46}
\begin {aligned}
&\quad\min(\overline{\al}_i, \overline{\al}_j) \sum_{j\in\E_i}\frac{1}{h_i}\|b_{ij}\|_{L^2(\pa\om_i)}^2\\
&\leq C\min(\overline{\al}_i, \overline{\al}_j) \sum_{j\in\E_i}\frac{1}{h_i}\|(w_j - \hat{w}_j)_j - (w_j - \hat{w}_j)_i\|_{L^2(E_{ij})}^2\\
&\leq C\Lr{\frac{H}{h}}^2\sum_{j\in\E_i}\frac{\overline{\al}_j} {\underline{\al}_j}d_j(\H_j w_j, \H_j w_j),
\end {aligned}
\end {equation}
since $H_i\simeq H_j$. Combining the inequalities~\eqref {app2esta22},~\eqref {app2esta43},~\eqref {app2esta45} and~\eqref {app2esta46}, we have
\begin {equation}\label {app2estofI1}
\begin {aligned}
I_1\leq C \Lr{\frac{H}{h}}^2\max_{i=1}^N \frac{\overline{\al}_i} {\underline{\al}_i} \sum_{j\in\E_i}\sum_{k=i,j}d_k(\H_k w_k, \H_k w_k).
\end {aligned}
\end {equation}
Substituting~\eqref {app2estofI1} and~\eqref {estimateofI2} into~\eqref {estimate}, we get
\begin {equation}\label {esapproach2}
\begin {aligned}
d_i(\H_i v_i, \H_i v_i)\leq C \Lr{\frac{H}{h}}^2\max_{i=1}^N \frac{\overline{\al}_i} {\underline{\al}_i} \sum_{j\in\E_i}\sum_{k=i,j}\overline{\al}_k d_k(\H_k w_k, \H_k w_k).
\end {aligned}
\end {equation}
By the summation of the above inequality for all $1\leq i \leq N$ and noting that the number of edges of each subdomain can be bounded independently of $N$, we finally obtain~\eqref {normofP} with $\beta$ satisfying~\eqref {beta1}.

Next we consider the special case when the coefficient $\al(x)$ in the subdomains $\om_i$ satisfies~\eqref {coefcond} for all $1\leq i\leq N$.
\begin {equation}\label {I1step2}
\begin {aligned}
I_1 &\leq C \min(\overline{\al}_i, \overline{\al}_j)\{|c_{ij}|_{H^{1/2}(\pa\om_i)}^2 + |d_{ij}|_{H^{1/2}(\pa\om_i)}^2+\\
&\quad+ \frac{1}{h_i}\|c_{ij}\|_{L^2(\pa\om_i)}^2  + \frac{1}{h_i}\|d_{ij}\|_{L^2(\pa\om_i)}^2\},
\end {aligned}
\end {equation}
where
\begin {align*}
c_{ij}=(w_i-\hat{w}_i)_i - (w_j-\hat{w}_j)_j\quad\text{and}\quad d_{ij}=
(w_j-\hat{w}_j)_j - (w_j-\hat{w}_j)_i.
\end {align*}
It is well-known that; c.f.~\cite {TW05},
\begin {equation}\label {esta1}
\begin {aligned}
&\quad\min(\overline{\al}_i, \overline{\al}_j)|c_{ij}|_{H^{1/2}(\pa\om_i)}^2\\
&\leq\min(\overline{\al}_i, \overline{\al}_j)\Lr{|(w_i-\hat{w}_i)_i|_{H^{1/2}(\pa\om_i)}^2
+ |(w_j-\hat{w}_j)_j|_{H^{1/2}(\pa\om_i)}^2}\\
&\leq\min(\overline{\al}_i, \overline{\al}_j)\sum_{j\in\E_i}\Lr{|(w_i-\hat{w}_i)_i|_{H^{1/2}(E_{ij})}^2
+ |(w_j-\hat{w}_j)_j|_{H^{1/2}(E_{ij})}^2}\\
&\leq C\sum_{j\in\E_i}\sum_{k=i,j}\overline{\al}_k(1+\log\frac{H_k}{h_k})^2\|\na \accentset{\circ}{\H}_k(w_k)_k\|_{L^{2}(\om_k)}^2\\
&= C\sum_{j\in\E_i}\sum_{k=i,j}\overline{\al}_k(1+\log\frac{H_k}{h_k})^2\|\na \H_k w_k\|_{L^{2}(\om_k)}^2\\
&\leq C \sum_{j\in\E_i}\sum_{k=i,j}\frac{\overline{\al}_k}{\underline{\al}_k}(1+\log\frac{H_k}{h_k})^2 d_k(\H_k w_k, \H_k w_k),
\end {aligned}
\end {equation}
since $\al(x)\geq \underline{\al}_k$ for all $x\in\om_k$.

Using (4.44) in~\cite {DGS13}, we have
\begin {equation}\label {esta2}
\begin {aligned}
&\quad\min(\overline{\al}_i, \overline{\al}_j)|d_{ij}|_{H^{1/2}(\pa\om_i)}^2\\
&=\min(\overline{\al}_i, \overline{\al}_j)\sum_{j\in\E_i}|(w_j-\hat{w}_j)_j - (w_j-\hat{w}_j)_i|_{H^{1/2}(E_{ij})}^2\\
&\leq C(1+\log\frac{H}{h})^2\sum_{j\in\E_i}\Lr{\overline{\al}_j\|\na\H_j w_j\|_{L^2(\om_j)}^2
+ \max(\frac{\overline{\al}_i}{\underline{\al}_i}, \frac{\overline{\al}_j}{\underline{\al}_j})\frac{1}{h_{ij}}\|\al_{ij}^{1/2}[(w_j)_i-(w_j)_j]\|_{L^2(E_{ji})}^2}\\
&\leq C(1+\log\frac{H}{h})^2\max_{i=1}^N \frac{\overline{\al}_i}{\underline{\al}_i}\sum_{j\in\E_i}d_j(\H_j w_j, \H_j w_j),
\end {aligned}
\end {equation}
where we have used~\eqref {ineqaver}, and the fact that $\delta$ is practically chosen such that $\delta=O(1)$.

Proceeding with the same lines of~\eqref {esta1}, we can obtain
\begin {equation}\label {esta3}
\begin {aligned}
&\quad\min(\overline{\al}_i, \overline{\al}_j)\frac{1}{h_i}\|c_{ij}\|_{L^2(\pa\om_i)}^2\\
&\leq\min(\overline{\al}_i, \overline{\al}_j)\frac{1}{h_i}\sum_{j\in\E_i}\Lr{\|(w_i-\hat{w}_i)_i\|_{L^2(E_{ij})}^2
+ \|(w_j-\hat{w}_j)_j\|_{L^2(E_{ij})}^2}\\
&\leq\min(\overline{\al}_i, \overline{\al}_j)\frac{H_i}{h_i}\sum_{j\in\E_i}\Lr{|(w_i-\hat{w}_i)_i|_{H^{1/2}(E_{ij})}^2
+ |(w_j-\hat{w}_j)_j|_{H^{1/2}(E_{ij})}^2}\\
&\leq C \frac{H}{h}(1+\log\frac{H}{h})^2\sum_{j\in\E_i}\sum_{k=i,j}\frac{\overline{\al}_k}{\underline{\al}_k} d_k(\H_k w_k, \H_k w_k),
\end {aligned}
\end {equation}
since $(w_i-\hat{w}_i)_i = 0$ at the end points of $E_{ij}$.

Using the inverse inequality, and the $L_2$ stability of the $L_2$ projection we have
\begin {equation}\label {esta41temp}
\begin {aligned}
&\quad\|(w_j-\hat{w}_j)_j - (w_j-\hat{w}_j)_i\|^2_{L^2(E_{ij})}\\
&\leq C\{\|Q_i[(w_j)_j - (w_j)_i]\|^2_{L^2(E_{ij})} + \|Q_i(w_j-\hat{w}_j)_j\|^2_{L^2(E_{ij})}+\\
&\quad+\|(w_j-\hat{w}_j)_j\|^2_{L^2(E_{ij})}+\|(\hat{w}_j)_i - (\hat{w}_j)_j\|^2_{L^2(E_{ij})}\}\\
&\leq C\{\|(w_j)_j - (w_j)_i\|^2_{L^2(E_{ij})}+\|(w_j-\hat{w}_j)_j\|^2_{L^2(E_{ij})}+\|(\hat{w}_j)_i - (\hat{w}_j)_j\|^2_{L^2(E_{ij})}\}\\
&\leq C\{\|(w_j)_j - (w_j)_i\|^2_{L^2(E_{ij})} + H_i |(w_j-\hat{w}_j)_j|^2_{H^{1/2}(E_{ij})} + H_i \max_{\partial E_{ij}}((w_j)_i - (w_j)_j)^2\}\\
&\leq C\{\|(w_j)_j - (w_j)_i\|^2_{L^2(E_{ij})} + H_i (1+\log \frac{H_j}{h_j})^2\|\nabla\H_j w_j\|^2_{L^2(\om_j)}+\\
&\quad+\frac{H_i}{h_i}\|(w_j)_j - (w_j)_i\|^2_{L^2(E_{ij})} + H_i(1+\log\frac{H}{h})\|\nabla\H_j w_j\|^2_{L^2(\om_j)}\},
\end {aligned}
\end {equation}
where we have used (4.43) in~\cite {DGS13}. Hence,
\begin {equation}\label {esta41}
\begin {aligned}
&\quad\min(\overline{\al}_i, \overline{\al}_j)\frac{1}{h_i}\|(w_j-\hat{w}_j)_j - (w_j - \hat{w}_j)_i\|_{L^2(E_{ij})}^2\\
&\leq C \{\frac{H_i}{h_i}\max(\frac{\overline{\al}_i}{\underline{\al}_i}, \frac{\overline{\al}_j}{\underline{\al}_j})\frac{1}{h_{ij}}\|\al_{ij}^{1/2}[(w_j)_j - (w_j)_i]\|^2_{L^2(E_{ij})} + \frac{H_i}{h_i}(1+\log\frac{H}{h})^2\frac{\overline{\al}_j}{\underline{\al}_j}\|\na\H_j w_j\|^2_{L^2(\om_j)}\}.
\end {aligned}
\end {equation}
This immediately gives
\begin {equation}\label {esta4}
\begin {aligned}
\min(\overline{\al}_i, \overline{\al}_j)\frac{1}{h_i}\|d_{ij}\|_{L^2(\pa\om_i)}^2\leq C \frac{H}{h}(1+\log\frac{H}{h})^2\max_{i=1}^N\frac{\overline{\al}_i}{\underline{\al}_i}\sum_{j\in\E_i} d_j(\H_j w_j, \H_j w_j).
\end {aligned}
\end {equation}
With the same arguments as in~\eqref {app2estofI1} and~\eqref {esapproach2}, we finally obtain~\eqref {normofP} with $\beta$ satisfying~\eqref {beta2}.

\end {proof}

\section {Numerical Experiments}

Let the domain $\om$ be a unit square $(0,1)^2$. For the experiments, we partition the domain $\om$ into $4\times 4$ square subdomains. The distribution of coefficients in each example is presented by figures. We use the proposed FETI-DP method for the discontinuous Galerkin formulation (Section \ref{FETI-DP for DG}) of the problem, and iterate with the preconditioned conjugate gradient (PCG) method. The iteration in each test stops whenever the $l_2$ norm of the residual is reduced by a factor of $10^{-6}$. The penalty parameter $\delta$ is chosen to be $5$ in all the experiments.

\begin {figure}[htbp]
\vskip -1.5in
\centering
\begin{minipage}[c]{0.49\textwidth}
\centering
\includegraphics[width=3.6in]{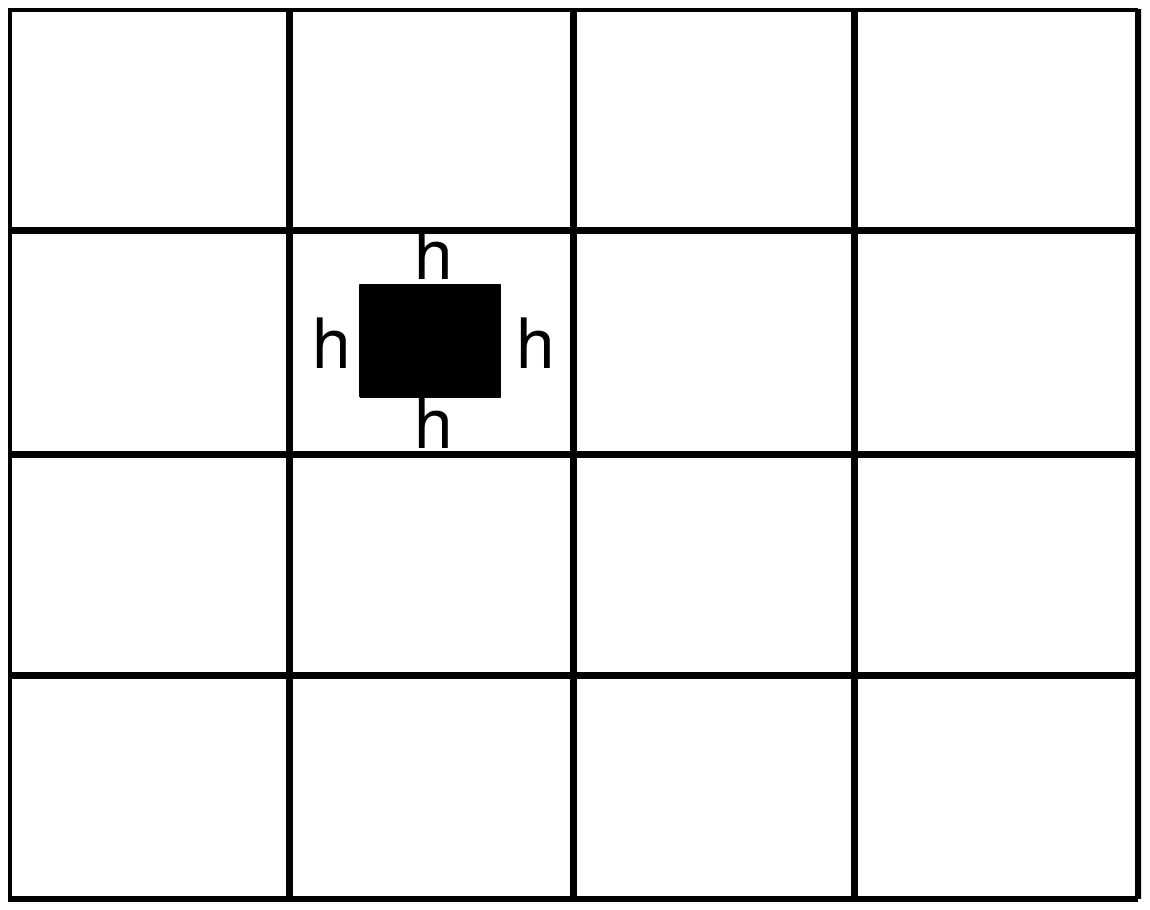}
\end{minipage}
\begin{minipage}[c]{0.49\textwidth}
\centering
\includegraphics[width=3.6in]{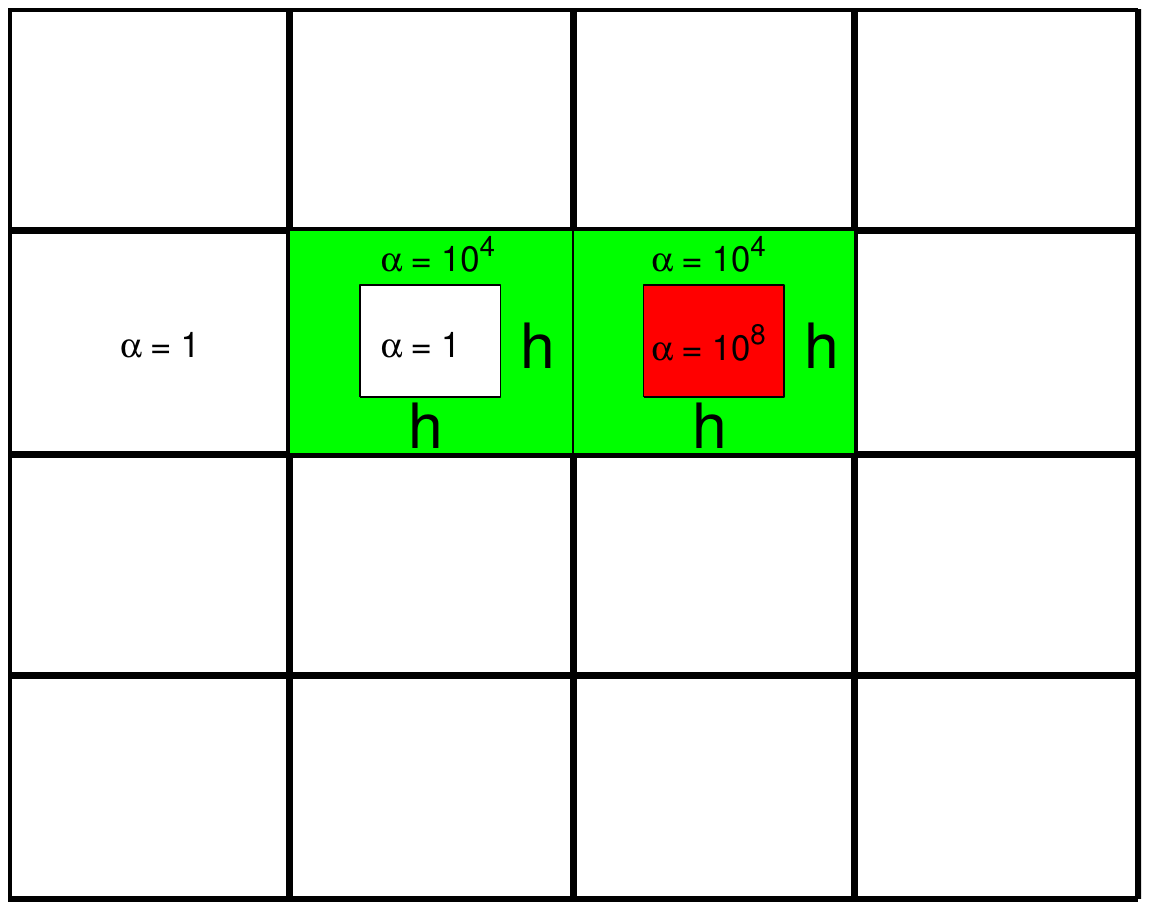}
\end{minipage}
\vskip -1.5in

\caption{\small Subdomain partition and coefficient distribution. Left: Example~\ref {example1}; Right: Example~\ref {example2}.}\label {figex1}
\end {figure}

\begin {example}\label {example1} In our first example, c.f. left picture of Fig.~\ref {figex1}, the coefficient denotes a 'binary' medium with $\alpha(x) = \widehat{\alpha}$ on a square shaped inclusion (shaded region) lying inside one subdomain $\om_i$ at a distance of $h$ from both the horizontal and the vertical edges of $\partial\om_i$, and $\alpha(x)=1$ in the rest of the domain. We study the behavior of the preconditioner as $h$ and $\widehat{\alpha}$ varies, respectively.
\end {example}

It follows from Tab.~\ref {tabex1} that the condition numbers are independent of the values of $\widehat{\alpha}$ since the coefficient contrast in the boundary layer is exactly equal to $1$. This is consistent with our theoretical results.

\begin {figure}[htbp]
\vskip -1.5in
\centering
\begin{minipage}[c]{0.49\textwidth}
\centering
\includegraphics[width=3.6in]{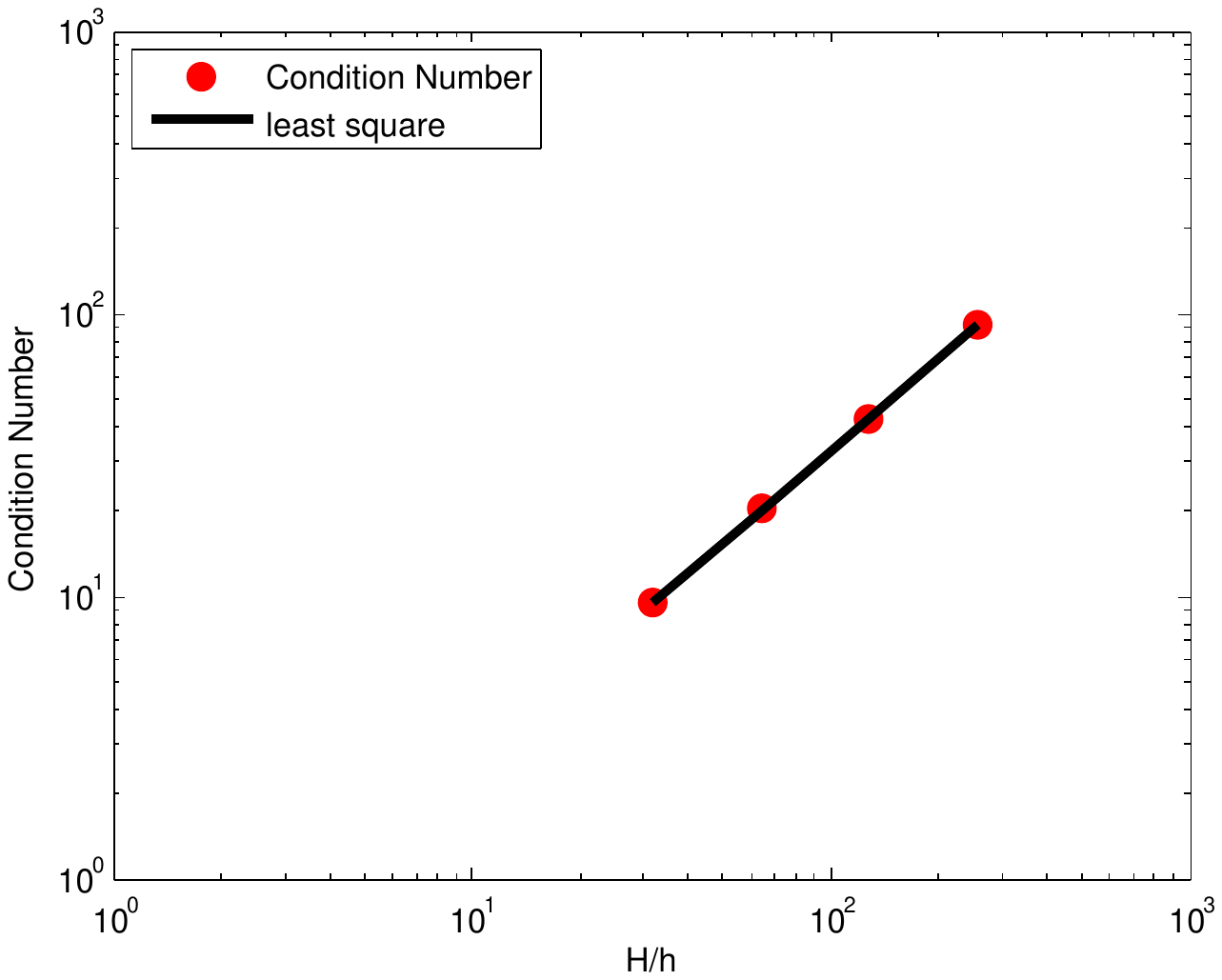}
\end{minipage}
\begin{minipage}[c]{0.49\textwidth}
\centering
\includegraphics[width=3.6in]{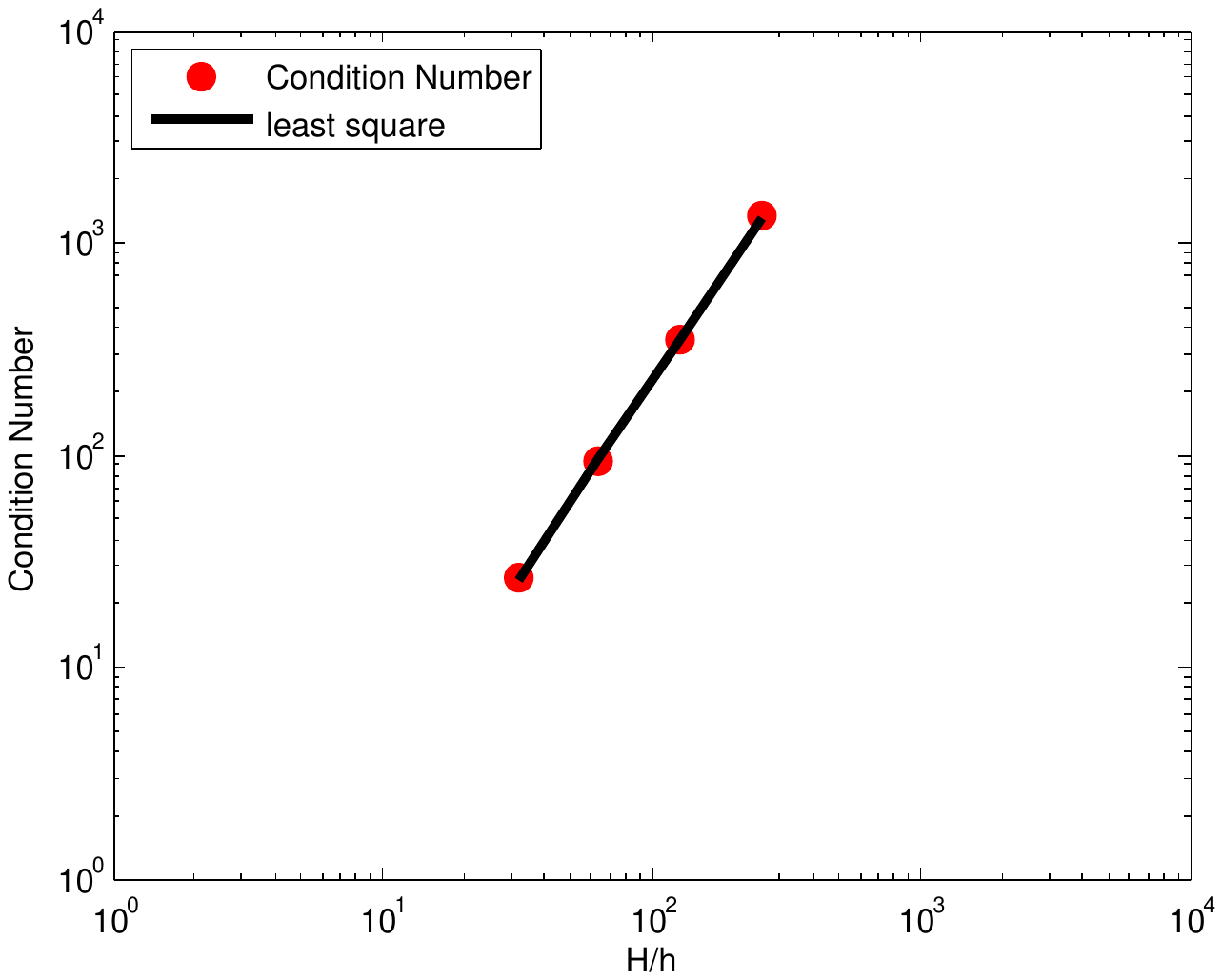}
\end{minipage}
\vskip -1.5in

\caption{\small Log-log plot of condition numbers vs. $H/h$. Left: Example~\ref {example1} with $\widehat{\alpha} = 10^6$, the slope of least square is $1.1$; Right: Example~\ref {example2}, the slope of least squre is $1.9$.}\label {figlog}
\end {figure}

Adopting different fine mesh sizes $h$, we obtain the log-log plot of the condition numbers in terms of $H/h$ for $\widehat{\alpha} = 10^6$. The left plot of Fig.~\ref {figlog} shows a dependence worse than linear growth, which is expected to become harder as $h$ goes finer. This confirms the estimate of~\eqref {beta1} that contains a logarithmic factor besides the linear dependence.

\begin {example}\label {example2} The distribution of coefficient is shown in the right picture of Fig.~\ref {figex1}, with inclusions in two neighbouring subdomains with coefficient values both larger and smaller than in the boundary layers.
\end {example}

Similar as the above example, we investigate the dependence of the condition numbers on the mesh ratio $H/h$. The right plot of Fig.~\ref {figlog} tells us the robustness of the quadratic dependence in the estimate of~\eqref {beta2}.

\begin {table}[htbp]
\centering

\caption{Example~\ref {example1}: PCG iterations and condition numbers (in parentheses).} \label {tabex1}
\begin {tabular}{cccccccc}
  \hline
  \noalign{\smallskip}

 &$H=32h$ & $H=64h$ & $H=128h$ & $H=256h$   \\
\hline
$\widehat{\alpha} = 10^2$& 13(8.568) &18(17.39) & 22(31.91)& 27(55.93)\\
$\widehat{\alpha} = 10^4$& 13(9.470) &17(20.30) & 22(42.39)& 29(89.93) \\
$\widehat{\alpha} = 10^6$& 13(9.481) &19(20.34) & 22(42.58)& 29(90.72)\\
\hline
\end {tabular}
\end {table}

\begin {table}[htbp]
\centering

\caption{Example~\ref {example2}: PCG iterations and condition numbers (in parentheses).} \label {tabex2}
\begin {tabular}{cccccccc}
  \hline

 $H=32h$ & $H=64h$ & $H=128h$ & $H=256h$   \\
\hline
 19(26.51) &25(92.54)  & 36(346.3)& 57(1333)\\
\hline
\end {tabular}
\end {table}

\begin {figure}[htbp]
\vskip -1.5in
\centering
\includegraphics[width=3.6in]{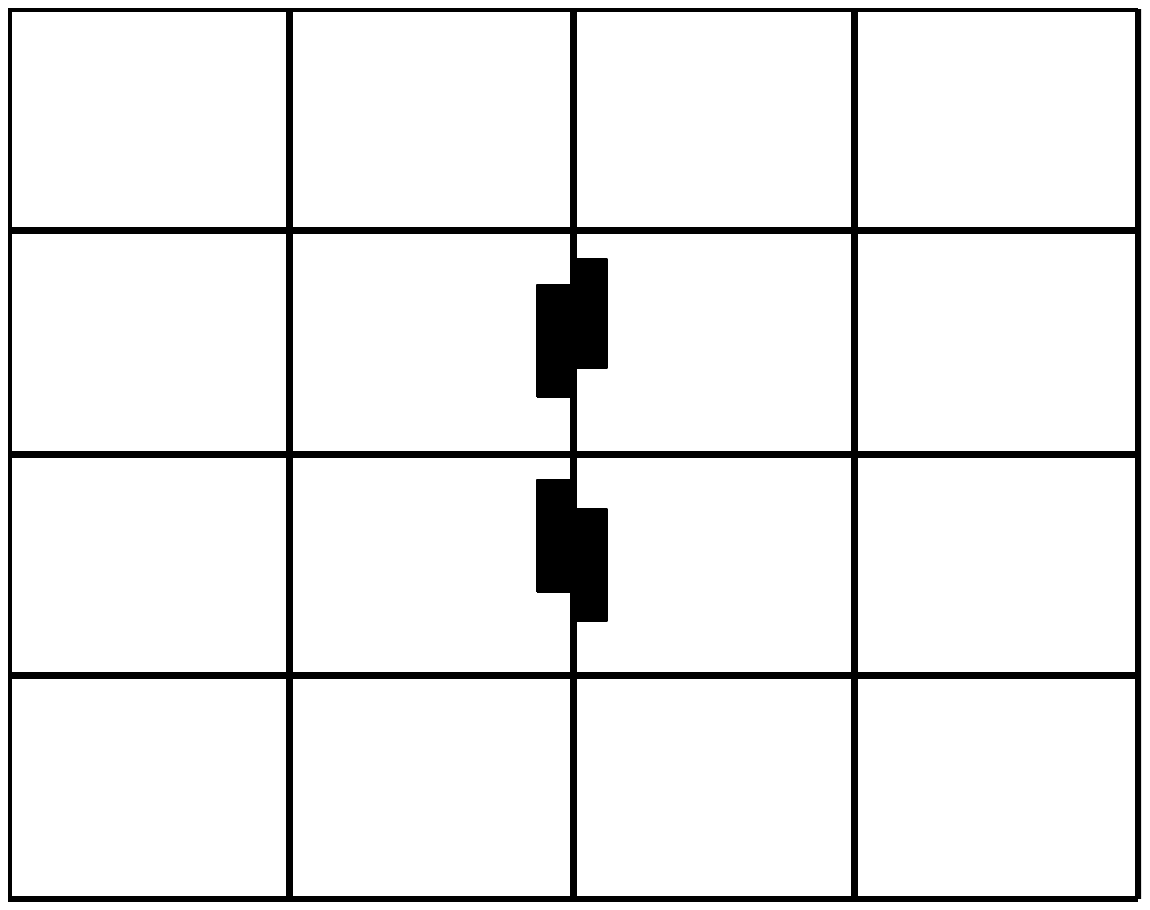}
\vskip -1.5in

\caption{\small Example~\ref {example3}: subdomain partition and coefficient distribution. The length of each inclusion is $H/8$ and the height is $H/2$.}\label {figex3}
\end {figure}

\begin {example}\label {example3} We employ this example to investigate the dependence of our method on the coefficient variation in the boundary layers. The distribution of the coefficient is depicted in Fig.~\ref {figex3}. The coefficient $\alpha(x) = \widehat{\alpha}$ in the edge islands (shaded region), and $\alpha(x) = 1$ else where.
\end {example}

The numerical results reported in Tab.~\ref {tabex3} confirm our theoretical results in Theorem~\ref {mainthm}, i.e., a linear dependence of the condition number on the coefficient variation in the boundary layers. It is worth further investigation to provide techniques to remove this dependence. In~\cite {PS08}, the authors used a pointwise weight to define the scaling matrix and finally made the performance of the method completely independent of the coefficient contrast for some special cases. However, there was no theoretical support to explain this robustness and this technique is not valid for the present example either.

%\begin {table}[htbp]
%\centering
%\renewcommand{\tablename}{Tab.}
%\caption{Example~\ref {example3}: PCG iterations and condition numbers (in parentheses).} \label {tabex3}
%\begin {tabular}{ccccccccc}
%  \hline
%$\widehat{\alpha}$ & $10^0$ & $10^1$ & $10^2$ & $10^3$&$10^4$&$10^5$&$10^6$\\
%\hline
% $H=64h$ &11(3.999)  & 22(11.95) &44(64.63) & 66(6.37$e+$2) & 91(6.36$e+$3) & 121(6.36$e+$4) & 145(6.36$e+$5) \\
%\hline
%\end {tabular}
%\end {table}

\begin {table}[htbp]
\centering

\caption{Example~\ref {example3}: PCG iterations and condition numbers (in parentheses).} \label {tabex3}
\begin {tabular}{ccccccccc}
  \hline
$\widehat{\alpha}$ & $10^2$ & $10^3$&$10^4$&$10^5$&$10^6$\\
\hline
 $H=64h$ &44(64.63) & 66(6.37$e+$2) & 91(6.36$e+$3) & 121(6.36$e+$4) & 145(6.36$e+$5) \\
\hline
\end {tabular}
\end {table}
%-------------------------REFERENCE-------------
\bibliographystyle{amsplain}

\providecommand{\bysame}{\leavevmode\hbox to3em{\hrulefill}\thinspace}
\providecommand{\MR}{\relax\ifhmode\unskip\space\fi MR }
% \MRhref is called by the amsart/book/proc definition of \MR.
\providecommand{\MRhref}[2]{%
  \href{http://www.ams.org/mathscinet-getitem?mr=#1}{#2}
}
\providecommand{\href}[2]{#2}
\begin{thebibliography}{}

\end{thebibliography}


\begin{thebibliography}{10}
\bibitem{Af03} R.A. Adams and J.J.F. Fournier, {\em Sobolev Spaces}, 2nd edition,
Academic Press, New York, 2003.

\bibitem{ABCM02} D.N. Arnold, F. Brezzi, B. Cockburn, and L.D. Marini, {\em Unified analysis of discontinuous Galerkin methods for elliptic problems}, SIAM J. Numer. Anal. \textbf{39}(2002), 1749--1779.

\bibitem{BS02} Susanne C. Brenner and L. Ridgway Scott,
\emph{The mathematical theory of finite element methods}, Springer, 2002.

\bibitem{Do03} C.R. Dohrmann, {\em A preconditioner for substructuring based on constrained energy minimization}, SIAM J. Sci. Comput. \textbf{25}(2003), 246--258.

\bibitem{DDG03} M. Dryja, {\em On discontinuous Galerkin methods for elliptic problems with discontinuous coefficients}, Comput. Methods Appl. Math.
\textbf{3}(2003), 76--85.

\bibitem{DGS07} M. Dryja, J. Galvis, and M. Sarkis, {\em BDDC methods for discontinuous Galerkin discretization of elliptic problems}, J. Complexity \textbf{23}(2007), 715--739.

\bibitem{DGS13} M. Dryja, J. Galvis, and M. Sarkis, {\em A FETI-DP preconditioner for a composite finite element and discontinuous Galerkin method}, SIAM Journal of Numerical Analysis \textbf{51}(2013), 400--422.

\bibitem{E92} W. E., {\em Homogenization of linear and nonlinear transport equations}, Comm. Pure Appl. Math. \textbf{XLV}(1992), 301--326.

\bibitem{ED02} Y. Efendiev and L.J. Durlofsky, {\em Numerical modelling of subgrid heterogeneity in two phase flow simulations}, Water resour. Res. \textbf{38}(2002), 1128--1138.

 \bibitem{EH09} Y. Efendiev and T. Hou, {\em Multiscale finite element methods}, Springer, 2009.

\bibitem{EGLMS13} Y. Efendiev, J. Galvis, R. Lazarov, M. Moon, and M. Sarkis, {\em Generalized Multiscale Finite Element Method. Symmetric Interior Penalty Coupling
}, preprint.

\bibitem{FLP00} C. Farhat, M. Lesoinne, and K. Pierson, {\em A scalable dual-primal domain decomposition method}, Numer. Linear Algebra \textbf{7}(2000), 687--714. Preconditioning techniques for large sparse matrix problems in industrial applications, Minneapolis, MN, 1999.

\bibitem{FL01} C. Farhat, M. Lesoinne, P. Le Tallec, K. Pierson, and D. Rixen, {\em FETI-DP: A dual primal unified FETI method I: a faster alternative to the two-level FETI method}, Int. J. Numer. Methods Eng. \textbf{50}(2001), 1523--1544.

\bibitem{FMR94} C. Farhat, J. Mandel, and F.X. Roux, {\em Optimal convergence properties of the FETI domain decomposition method}, Comput. Meth. Appl. Mech. Engrg. \textbf{}(1994), 367--388.

\bibitem{FR91} C. Farhat and F.X. Roux, {\em A method of finite element tearing and interconnecting and its parallel solution algorithm}, Int. J. Numer. Methods Eng. \textbf{}(1991), 1205--1277.


\bibitem{HP02} B. Heinrich and K. Pietsch, {\em Nitsche type mortaring for some elliptic problem with corner singularities}, Computing \textbf{68}(2002), 217--238.

\bibitem{KWD02} A. Klawonn, O. Widlund, and M. Dryja, {\em Dual-primal FETI methods for three-dimensional elliptic problems with heterogenenous coefficients}, SIAM J. Numer. Anal. \textbf{40}(2002), 159--179.

\bibitem{LJ06} I. Lunati and P. Jenny, {\em Multiscale finite-volume method for compressible multiphase flow in porous media}, Journal of Computational Physics \textbf{216}(2006), 616--636.

\bibitem{LJ08} I. Lunati and P. Jenny, {\em Multiscale finite-volulme method for density-driven flow in porous media}, Comput. Geosci. \textbf{12}(2008), 337--350.

\bibitem{MD03} J. Mandel and C.R. Dohrmann, {\em Convergence of a balancing domain decomposition by constraints and energy minimization}, numer. Lin. Alg. Appl. \textbf{10}(2003), 639--659.

 \bibitem{MT01} J. Mandel and R. Tezaur, {\em On the convergence rate of a dual-primal substructuring method}, Numer. Math. \textbf{88}(2001), 543--558.

\bibitem{PS08} C. Pechstein and R. Scheichl, {\em Analysis of FETI methods for multiscale PDEs}, Numer. Math. \textbf{111}(2008), 293--333.

\bibitem{PS11} C. Pechstein and R. Scheichl, {\em Analysis of FETI methods for multiscale PDEs. Part II: interface variation}, Numer. Math. \textbf{118}(2011), 485--529.

\bibitem{Ri08} B. Rivi$\grave{e}$re, {\em Discontinuous Galerkin methods for solving elliptic and parabolic equations: Theory and implementations}, vol. 35 of Frontiers in Applied Mathematics, SIAM, Philadelphia, PA, 2008.

\bibitem{St98} R. Stenberg, {\em Mortaring by a method of J.A.Nitsche}, in Computational mechanics (Buenos Aires, 1998), Centro Internac. M$\acute{e}$todos Num$\acute{e}$r. Ing., Barcelona, 1998, CD-ROM file.

\bibitem{SBG78}B.F. Smith, P. Bj{\o}rstad, and W. Gropp,
\emph{Domain Decomposition: Parallel Multilevel Methods for Elliptic Partial Differential Equations}, Cambridge University Press, 1996.

\bibitem{TW05}A. Toselli and O. Widlund,
\emph{Domain Decomopsition Methods-Algorithms and Theory}, Springer-Verlag Berlin Heidelberg, Berlin, 2005.









%
%
%
%
%\bibitem{T.X.Z.}Tarek P. A. Mathew,\,\, Domain decomposition methods for the numerical solution of partial differential equations, volume 61 of Lecture Notes in Computational Science and Engineering,\,\, Springer-Verlag,\,\, Berlin, \,\,2008.
%\bibitem{B.D.V.} Petter Bj{\o}rstad,\,\,Maksymilan Dryja,\,\,Eero Vainikko,\,\, Additive Schwarz Methods without Subdomain Overlap and New Coarse Spaces,\,\,Proceedings from the
%              8th. International conference on domain decomposition methods,\,\,1996,\,\,John Wiley \& Sons,\,\,1997.
%\bibitem{D.G.S.1}Maksymilian Dryja,\,\,Juan Galvis,\,\,Marcus Sarkis,\,\,BDDC methods for discontinuous Galerkin discretization of elliptic problems,\,\,J. Complexity,\,\,2007,\,\,715-739.
%\bibitem{D.S.1}M. Dryja,\,\, M. Sarkis,\,\, Additive average Schwarz methods for discretization of elliptic problems with highly discontinuous coefficients,\,\,
%Computational Methods in Applied Mathematics,\,\, 10(2010),\,\,1-13.
%\bibitem{D.S.2}Maksymilian Dryja,\,\,Marcus Sarkis,\,\,FETI-DP method for DG discretization of elliptic problems with discontinuous coefficients,\,\,2010, In Review.
%\bibitem{D.G.S.}M. Dryja,\,\, J. Galvis,\,\, M. Sarkis,\,\, Neumann-Neuamnn methods for a DG discretization on geometrically nonconforming substructures,\,\,To appear(2011).
%\bibitem{D.G.S.2}Maksymilian Dryja, Juan Galvis, Marcus Sarkis,\,\, Neumann-Neumann methods for a DG discretization of elliptic problems with discontinuous coefficients on geometrically nonconforming substructures,\,\, Numerical Methods for Partial Differential Equations,\,\, 2011.
%\bibitem{I.P.R.} I. G. Graham,\,\, P. O. Lechner,\,\, R. Scheichl,\,\, Domain Decomposition for multiscale PDEs,\,\, Numer. Math.,\,\, 2007,\,\, 589-626.
%\bibitem{C.G.S.S.}Cliffe, K. A.,\,\, Graham, I. G.,\,\, Scheichl, R.,\,\, Stals, L., \,\,Parallel computation of flow in heterogeneous media modelled by mixed finite elements.\,\, J. Comp. Phys.,\,\, 164(2000),\,\, 258-282.
%
%\bibitem{T.X.Z.b} Thomas Y. Hou,\,\,Xiaohui Wu,\,\,Zhiqiang Cai,\,\, Convergence of a multiscale finite element method for elliptic problems with rapidly oscillating coefficients,\,\,1999,\,\,913-943.
%\bibitem{G.E.}J. Galvis, Y. Efendiev, \,\,Domain decomposition preconditioners for multiscale flows in high-contrast media,
%              \,\,Multiscale Model Simul.,\,\,8(2010), \,\,1461-1483.
%
%\bibitem{P.S.1}Clements Pechstein, \,\,Robert Scheichl,\,\,Analysis of FETI methods for multiscale pdes,\,\,Numer. Math.,\,\,2008, 293-333.
%\bibitem{P.S.2}Clements Pechstein, \,\,Robert Scheichl,\,\,Analysis of FETI methods for multiscale pdes-Part $\pi$:interface variation,\,\,Numer. Math.,\,\,2011, 485-529.
%
%
%\bibitem{F.M.R.}Charbel Farhat,\,\, Jan Mandel,\,\,Francois Xavier Roux,\,\, Optimal convergence propertities of the FETI domain decomposition method,\,\,Comput. Methods Appl. Mech. Engrg.,\,\,1994,\,\,367-388.
%\bibitem{I.P.} Ivan Lunati,\,\,Patrick Jenny,\,\,Multiscale finite-volulme method for density-driven flow in porous media,\,\,Comput. Geosci.,\,\,2008,\,\,337-350.
%\bibitem{I.P.b} Ivan Lunati,\,\,Patrick Jenny,\,\,Multiscale finite-volume method for compressible multiphase flow in porous media,\,\,Journal of Computational Physics,\,\,2006,\,\,616-636.
%\bibitem{F.L.T.P.R.}Charbel Farhat,\,\,Michel Lesoinne,\,\,Patrick Le Tallec,\,\,Kendall Pierson,\,\,Daniel Rixen,\,\,FETI-DP:A dual primal unified FETI method-part I:A faster alternative to the two level FETI method,\,\,Internat. J. Numer. Methods Engrg.,\,\,2001,\,\,1523-1544.
%\bibitem{F.R.}Charbel Farhat,\,\, Francois Xavier Roux,\,\, A method of finite element tearing and interconnecting and its parallel solution algorithm,\,\,Internat. J. Numer. Methods Engrg.,\,\,1991,\,\,1205-1227.
%\bibitem{Z.X.S.}Zhiqiang Cai,\,\, Xiu Ye,\,\, Shun Zhang,\,\, Discontinuous Galerkin finite element methods for interface
%problems: a priori and a posteriori error estimations,\, SIAM J. Numer. Anal.,\,\,49(2011).
%\bibitem{X.O.} Xiaobing Feng,\,\, Ohannes A. Karakasshian,\,\, Two-level additive Schwarz methods for a discontinous Galerkin approximation of second order elliptic problems,\,\, SIAM J. Numer. Anal.,\,\,(39)2001,\,\,1343-1365.
%
%\bibitem{M.T.}Jan Mandel,\,\,Radek Tezaur,\,\,On the convergence of a dual-primal substructuring method,\,\,Numer. Math.,\,\,88(2001),\,\,543-558.
%
%\bibitem{Gelhar}L.~W.~ Gelhar, A Stochastic Conceptual Analysis of One-Dimensional Groundwater
%        Flow in Nonuniform Homogeneous Media, Water Resour. Res.,\,\, {\bf 11},\,\, 725-741(1975).
%\bibitem{H.K.}R.~J. ~Hoeksema and P.~K.~ Kitanidis, Analysis of the Spatial Structure of Properties
%of Selected Aquifers, Water Resour. Res.,\,\, {\bf 21},\,\,  563-572 (1985).
%\bibitem{FR.DM.}Feng Ruan, and Dennis Maclaughlin, An efficient multivariate random field generator using the fast Fourier transform,\,\,Advances in Water Resources,\,\,{\bf 21},\,\, 385-399(1998).
%\bibitem{R.G.S.W.}M.~L.~ Robin, A.~L.~ Gutjahr, E.~A.~ Sudicky, and J.~L.~ Wilson, Cross-correlated Random
%Field Generation with the Direct Fourier Transform Method, Water Resour.
%Res.,\,\, {\bf 29} ,\,\, 2385-2398(1993).
%\bibitem{K.K.}Kozintsev, B., Kedem, B.: Gaussian package, University of Maryland. Available at http://www.
%math.umd.edu/ bnk/bak/generate.cgi (1999)~~\\
\end{thebibliography}

\end {document}